\documentclass[11pt]{article}

\usepackage{amsmath}
\usepackage{amsthm}
\usepackage{amssymb}
\usepackage{array}
\usepackage[ruled,vlined]{algorithm2e}
\usepackage[nohead,margin=1.1in]{geometry}
\usepackage{color}
\usepackage{verbatim}
\usepackage{booktabs}
\usepackage{threeparttable}
\usepackage{url}

\newtheorem{theorem}{Theorem}[section]

\newtheorem{proposition}{Proposition}[section]
\newtheorem{lemma}{Lemma}[section]
\newtheorem{corollary}{Corollary}[section]
\newtheorem{assumption}{Assumption}[section]

\newtheorem{remark}{Remark}[section]

\usepackage{graphicx,epstopdf}
\usepackage{subfigure}
\usepackage{subfig}
\graphicspath{{./pics/},{./result_v5/}}
\usepackage[toc,page]{appendix}
\allowdisplaybreaks

\numberwithin{equation}{section}
\newcommand{\E}{\mathbb{E}}

\begin{document}
\numberwithin{equation}{section}
\numberwithin{figure}{section}
\title{An Analysis of Stochastic Variance Reduced Gradient for Linear Inverse Problems}

\author{Bangti Jin\thanks{Department of Computer Science, University College London, Gower Street, London WC1E 6BT, UK
(b.jin@ucl.ac.uk,bangti.jin@gmail.com). The work of BJ is supported by UK EPSRC grant EP/T000864/1.}
\and Zehui Zhou\thanks{Department of Mathematics, The Chinese University of Hong Kong, Shatin, New Territories, Hong Kong.
({zhzhou@math.cuhk.edu.hk, zou@math.cuhk.edu.hk}).
The work of JZ was substantially supported by Hong Kong RGC General Research Fund (projects 14306718 and
14304517).} \and Jun Zou\footnotemark[2]}
\date{}
\maketitle

\begin{abstract}
Stochastic variance reduced gradient (SVRG) is a popular variance reduction technique for accelerating stochastic gradient descent (SGD).
We provide a first analysis of the method  for solving a class of linear inverse problems in the lens of the classical regularization
theory. We prove that for a suitable constant step size schedule, the method can achieve an optimal convergence rate in
terms of the noise level (under suitable regularity condition) and the variance of the SVRG iterate error is smaller than
that by SGD. These theoretical findings are corroborated by a set of numerical experiments.

\end{abstract}

\section{Introduction}

In this paper, we consider the numerical solution of the following finite-dimensional linear inverse problem:
\begin{equation}\label{eqn:lininv}
Ax=y^\dag,
\end{equation}
where $A\in \mathbb{R}^{n\times m}$ is the system matrix representing the data formation mechanism,
and $x\in \mathbb{R}^m$ is the unknown signal of interest. In practice, we only have access to a noisy
version $y^\delta$ of the exact data $y^\dag = Ax^\dag$ (with $x^\dag$ being the minimum norm solution
relative to the initial guess $x_0$, cf. \eqref{eqn:min-norm}), i.e.,
\begin{equation*}
y^\delta =y^\dag +\xi,
\end{equation*}
where $\xi \in \mathbb{R}^n$ denotes the noise in the data with a noise level $\delta =\|\xi\|$, with $\|\cdot\|$
being the Euclidean norm of a vector (and also the spectral norm of a matrix).
We denote the $i${th} row of the matrix $A$ by a column vector $a_i\in \mathbb{R}^m$, i.e., $A=[a_i^t]_{i=1}^n$
(with the superscript $t$ denoting the matrix/vector transpose), and the $i${th} entry of the vector
$y^\delta\in\mathbb{R}^n$ by $y_i^\delta $. Linear inverse problems of the form \eqref{eqn:lininv} arise in a
broad range of practical applications, e.g., computed tomography and optical imaging.

Over the last few years, stochastic iterative algorithms have received much interest in the inverse problems community.
The most prominent example is stochastic gradient descent (SGD) due to Robbins and Monro \cite{RobbinsMonro:1951}. The starting
point is the following optimization problem:
\begin{equation}\label{eqn:obj}
 J (x) = \tfrac{1}{2n}\|Ax-y^\delta\|^2 = \tfrac{1}{n}\sum_{i=1}^n f_i(x),\quad \mbox{with}\quad f_i(x)=\tfrac{1}{2} \big((a_i,x)-y_i^\delta\big)^2,
\end{equation}
where $(\cdot,\cdot)$ denotes the Euclidean inner product on $\mathbb{R}^m$.
Then SGD reads as follows. Given an initial guess
$\hat x_0^\delta \equiv x_0$, the iterate $\hat x_k^\delta$ is constructed as
\begin{equation*}
  \hat x_{k+1}^\delta = \hat x_k^\delta - \eta_k f'_{i_k}(\hat x_k^\delta),
\end{equation*}
where $\eta_k>0$ is the step size at the $(k+1)$-th step, and the index $ i_k$ is sampled uniformly from the index set $\{1,\ldots,n\}$.
One attractive feature of the method is that the computational complexity per iteration does not depend on the
data size $n$, and thus it is directly scalable to large data volume, which is especially attractive in the era
of big data. {SGD type methods have found applications in several inverse problems, e.g.,
randomized Kaczmarz method \cite{HermanLentLutz:1978,StrohmerVershynin:2009} in
computed tomography, ordered subset expectation maximization \cite{HudsonLarkin:1994,KeretaTwyman:2021}
for positron emission tomography, and more recently also some nonlinear inverse problems, e.g., optical
tomography \cite{ChenLiLiu:2018} and phonon transmission coefficient \cite{GambaLi:2020}.}

However, the relevant mathematical theory for inverse problems in the lens of regularization theory
\cite{EnglHankeNeubauer:1996,KaltenbacherNeubauerScherzer:2008,ItoJin:2015} is still not fully understood. Existing works \cite{JinLu:2019,JinZhouZou:2020,JahnJin:2020,JinZhouZou:2021} focus on the standard SGD for inverse problems,
proving that SGD is a regularization method when equipped with a suitable stopping criterion, and the SGD iterates
converge at a certain rate. However, the presence of stochastic gradient noise generally prevents SGD from
converging to the solution when a constant step size is used and leads to a slow, sublinear
rate of convergence when a diminishing {step size} schedule is employed.
Amongst various acceleration strategies, variance reduction (VR) represents one prominent idea that
has achieved great success, including SAG \cite{LeRouxSchmidt:2012}, SAGA \cite{DefazioBachLacoste:2014}, SVRG
\cite{JohnsonZhang:2013,ZhangMahdaviJin:2013} and SARAH \cite{NguyenLiu:2017} etc;
These methods take advantage of the finite-sum structure prevalent in machine learning problems, and exhibit
improved convergence behavior over SGD; see the work \cite{GowerSchmidtBach:2020}
for a recent overview of variance reduction techniques in machine learning.

Stochastic variance reduced gradient (SVRG) combines SGD with predictive variance reduction and is very popular
in stochastic optimization. It was proposed independently by two groups of researchers, i.e.,
Johnson and Zhang \cite{JohnsonZhang:2013} and Zhang, Mahdavi and Jin
\cite{ZhangMahdaviJin:2013}, for accelerating SGD for minimizing smooth and strongly convex objective functions.
When applied to problem \eqref{eqn:obj}, the basic version of SVRG reads as follows.
Given an initial guess $x_0^\delta \equiv x_0 \in \mathbb{R}^m$, SVRG updates the iterate $x_{k}^\delta $
recursively by
\begin{equation}\label{eqn:SGD}
  x_{k+1}^\delta =x_k^\delta -\eta_k \big( f'_{i_k}(x_k^\delta)- f'_{i_k}(x_{k_M}^\delta)+ J'(x_{k_M}^\delta)\big), \quad k=0,1,\cdots,
\end{equation}
where the row index $i_k$ is drawn uniformly from the index set $\{1,\cdots,n\}$, $\eta_k>0$ is the
step size at the $k$th iteration, $M$ is the frequency of computing the full gradient, and $k_M=
[\frac{k}{M}]M,$ ($[\cdot]$ takes the integral part of a real number). {The choice of
the frequency $M$ can affect the practical performance of the algorithm, and it was suggested to
be $2n$ and $5n$ for convex and nonconvex optimization, respectively \cite{JohnsonZhang:2013}. In this study, we
show that SVRG can achieve optimal convergence rates when $M$ is chosen such that $M\geq O(n^\frac12)$.} When compared with SGD in
\eqref{eqn:SGD}, SVRG employs the anchor / snapshot point $x_{k_M}^\delta$ to reduce the variance
of the gradient estimate: it computes the full gradient $J'(x_{k_M}^\delta)$ of $J$ at the anchor
point $x_{k_M}^\delta$ for every $M$ iterates, and then combines $J'(x_{k_M}^\delta)$ with
the gradient gap $f_{i_k}'(x_k^\delta)-f_{i_k}'(x_{k_M}^\delta)$ to obtain a new gradient estimate for
updating the SVRG iterate $x_{k+1}^\delta$. In contrast, SGD employs the stochastic gradient
$f_{i_k}'(\hat x_k^\delta)$ only, and the classical Landweber method uses only the gradient $J'(x)$. Thus,
SVRG can be viewed as a hybridization between the Landweber method and SGD. A detailed comparison between SGD
and SVRG are given in Algorithms \ref{alg:sgd} and \ref{alg:svrg}, where SVRG is stated in the
form of double loop. In practice, there are several variants of SVRG, dependent on the choice of
the anchor point, e.g., last iterate, iterate average, random choice and weighted iterate average
(within the inner loop). In this  work, we study only the version given in Algorithm \ref{alg:svrg}.

\begin{algorithm}[H]
\SetAlgoLined
Set initial guess $\hat x_0^\delta=x_0$ and step size {schedule $\eta_k$}\\
 \For{$k=0,1,\cdots$}{
draw $i_k$ i.i.d. uniformly from $\{1,\cdots,n\}$\\
update  $\hat x_{k+1}^\delta =\hat x_k^\delta -\eta_k ((a_{i_k} ,\hat x_k^\delta )-y_{i_k}^\delta )a_{i_k}$\\
check the stopping criterion\\
}
\caption{Stochastic Gradient Descent (SGD) for problem \eqref{eqn:lininv}.\label{alg:sgd}}
\end{algorithm}

\begin{algorithm}[H]
\SetAlgoLined
Set initial guess $x_0^\delta=x_0$, frequency $M$ and step size schedule $\eta_k$\\
 \For{$K=0,1,\cdots$}{
compute $ {J'(x_{KM}^\delta)}$\\
\For{$t=0,1,\cdots,M-1$}{
draw {$i_{KM+t}$} i.i.d. uniformly from $\{1,\cdots,n\}$\\
update $x_{KM+t+1}^\delta =x_{KM+t}^\delta -\eta_{KM+t} \big((a_{i_{KM+t}} ,x_{KM+t}^\delta-x_{KM}^\delta)a_{i_{KM+t}}+ {J'(x_{KM}^\delta)}\big)$\\
}
check the stopping criterion.
}
\caption{Stochastic Variance Reduced Gradient (SVRG) for problem \eqref{eqn:lininv}.\label{alg:svrg}}
\end{algorithm}

It is known that VR enables speeding up the convergence of the algorithm in the sense of optimization
\cite{BottouCurtisNocedal:2018,GowerSchmidtBach:2020}.
Since its first introduction, SVRG has received a lot of
attention within the optimization community, and several convergence results of SVRG and its variants
have been obtained \cite{HarikandehSchmidt:2015,AllenZhuHazan:2016,AllenZhuYuan:2016,ReddiHefny:2016,
XuLinYang:2017,KovalevRichtarik:2020,ShangJiao:2020}. Note that here the precise meaning of convergence
depends crucially on the property of the objective function $J(x)$: (i) the distance of the SVRG
iterate $x_k^\delta$ to a global minimizer for a strictly convex $J(x)$, (ii) the optimality
gap (i.e., $J(x_k^\delta)-\min_x J(x)$) for a convex $J(x)$ and (iii) the
norm of the gradient $\|J'(x_k^\delta)\|$ for a nonconvex $J(x)$, in terms of the iterate number $k$.
For example, Allen-Zhu and Hazan \cite{AllenZhuHazan:2016} proved that SVRG (with a different
choice of the anchor point) converges at an $O(n^\frac{2}{3}\epsilon^{-1})$ rate
to an approximate stationary point $x^*$ (i.e., $\|J'(x^*)\|^2\leq\epsilon$) for a nonconvex but smooth $J(x)$.
Reddi et al \cite{ReddiHefny:2016} proved a nonasymptotic rate of convergence of SVRG for nonconvex
optimization and identified a subclass of nonconvex problems (satisfied by gradient dominated functions)
for which a variant of SVRG attains linear convergence.

These important breakthroughs in the optimization literature naturally motivate the following question:
\textit{Does the desirable convergence property of SVRG carry over to inverse problems in the sense of
regularization theory?} {The answer to this question} is not self-evident, {since accelerated
iterative schemes do not necessarily retain the optimal convergence in the sense of regularization (see
\cite{Neubauer:2017,Kindermann:2021} for studies on Nesterov's accelerated scheme).}
For linear inverse problems in \eqref{eqn:lininv}, the objective $J(x)$ in \eqref{eqn:obj} is convex
but not strictly so. Further, it  is ill-posed in the sense that a global minimizer often does not
exist, and even if it does exists, it is unstable with respect to the inevitable perturbation of the data
$y^\delta$ and is probably physically irrelevant. Instead, we construct an approximate
minimizer that converges to the exact solution $x^\dag$ as the noise level $\delta$ tends to $0^+$ by
stopping the iteration properly, a procedure commonly known as iterative regularization  (by early stopping)
\cite{KaltenbacherNeubauerScherzer:2008}, and the accuracy of the approximation is measured in terms
of the noise level $\delta$. To the best of our knowledge, the theoretical properties of SVRG and
other variance reduction techniques have not been studied so far in the lens of regularization theory.

In this work, we contribute to the theoretical analysis of SVRG for a class of linear inverse problems from the
perspective of classical regularization theory \cite{EnglHankeNeubauer:1996,KaltenbacherNeubauerScherzer:2008,
ItoJin:2015}. Under the constant step size schedule and the canonical source condition, we prove that the epochwise SVRG iterate
$x_{KM}^\delta$ converges to the minimum norm solution $x^\dag$ at an optimal rate (in terms of $\delta$)
when combined with \textit{a priori} stopping rule, and that due to the built-in variance reduction mechanism,
for the same iterate number, the variance of SVRG iterate is indeed smaller than that of SGD, showing
the beneficial effect of variance reduction; see Theorems \ref{thm:main} and \ref{thm:svrg-sgd}. In particular,
SVRG allows using larger step sizes than that for SGD while still overcoming the undesirable saturation
phenomenon (cf. Remark \ref{rem:con_optimal}). See Section \ref{sec:main} for precise statements of the
theoretical findings and related discussions in the context of inverse problems. These theoretical results
are complemented by extensive numerical results in Section \ref{sec:numer}.

The rest of the paper is organized as follows. In Section \ref{sec:main} we present
and discuss the main results of the work. In Section \ref{sec:decomp}, we recall preliminary
results, especially a careful decomposition of the error of the epoch SVRG / SGD iterate into the bias and variance
components. In Section \ref{sec:conv} we give the convergence rate analysis, and
prove an optimal convergence rate, and in Section \ref{sec:comp} we present a comparative study of SVRG versus
SGD, and show that variance component of the SVRG error is smaller than that of the SGD error. Finally, in Section
\ref{sec:numer}, we present several numerical experiments to complement the theoretical analysis. {For better
readability, the lengthy and technical proofs of several auxiliary results are deferred to the appendix.} Throughout, the
notation $c$ with suitable subscripts denotes a generic constant.

\section{Main results and discussions}\label{sec:main}
In this section, we state the main results of the work. First we state the standing assumption.
We denote by $\mathcal{F}_k$ the filtration generated by the random indices $\{i_0,i_1,\ldots,
i_{k-1}\}$. Let $\mathcal{F}={\bigvee_{k=1}^\infty} \mathcal{F}_k$, $\mathcal{F}^c_k=\mathcal{F}
\setminus\mathcal{F}_k$, $(\Omega,\mathcal{F},\mathbb{P})$ being the associated probability space, and
$\E[\cdot]$ denotes taking the expectation with respect to the filtration $\mathcal{F}$ and $\E_{j}[\cdot]:=\E[\cdot|\mathcal{F}^c_{j+1}\cup\mathcal{F}_{j}]$. The
SVRG iterate $x_k^\delta$ is random, and measurable with respect to $\mathcal{F}_k$. Let $e_k^\delta
=x_k^\delta-x^\dag$ be the error of the SVRG iterate $x_k^\delta$ with respect to the unique
minimum-norm solution $x^\dag$, defined by
\begin{equation}\label{eqn:min-norm}
  x^\dag = \arg\min_{x\in\mathbb{R}^m: Ax = y^\dag} \|x-x_0\|.
\end{equation}
Let $B=\E[a_ia_i^t]=n^{-1}A^t A \in \mathbb{R}^{m\times m}$. Throughout we assume that $\|B\|\leq 1$,
which can easily be achieved by scaling. In this work we consider a constant step size
schedule, which is commonly employed by SVRG. Assumption \ref{ass:stepsize}(ii) is commonly known as
the source condition in the inverse problems literature \cite{EnglHankeNeubauer:1996}, which
implicitly assumes a certain regularity on the initial error. This condition is central for deriving
convergence rates. It is well known that in the absence of source type conditions, the convergence for
a regularization method can be arbitrarily
slow \cite{EnglHankeNeubauer:1996}. Assumption \ref{ass:stepsize}(iii) enables an important commuting
property (cf. Lemma \ref{lem:commut}), which greatly facilitates the analysis. Numerically
this property does not affect the performance of SVRG, and thus it seems largely due to the limitation
of the analysis technique.
\begin{assumption}\label{ass:stepsize}
The following assumptions hold.
\begin{itemize}
\item[$\rm(i)$] The step size $\eta_j = c_0$, $j=0,1,\cdots$, with $c_0\leq \big(\max(\max_i \|a_i\|^2,\|B\|^2)\big)^{-1}$.
\item[$\rm(ii)$] There exist some $\nu>0$ and $w\in\mathbb{R}^m$ such that the exact solution $x^\dag$ satisfies
$x^\dag-x_0=B^\nu w.$
\item[{\rm(iii)}] The matrix $A=\Sigma V^t$ with $\Sigma$ being diagonal and nonnegative and $V$ column orthonormal.
\end{itemize}
\end{assumption}

The next result represents the main theoretical contribution of the work. It implies
that SVRG can achieve the optimal convergence rate for linear inverse problems under
the given assumption on the step size. The step size restriction originates from the
fact that SVRG still employs a randomized gradient estimate for the iterate update,
albeit with reduced variance, when compared with the Landweber method. {Nonetheless,
the restriction on the step size is more benign than that for SGD: It allows
achieving optimal convergence rate under larger step size than that in SGD.}
\begin{theorem}\label{thm:main}
Let Assumption \ref{ass:stepsize} hold, and $c_*>1$ satisfy
\begin{align}\label{cond:sat_svrg}
&(4+2(Mc_0\|B\|)^2)nM^{-2}c_Bc_{B,M}\leq 1-c_*^{-1}\\
\mbox{with}\quad
&c_{B,M}=\sum_{i=1}^{M-1}(1-(1-c_0\|B\|)^i)^2 \mbox{ and } c_B=(1-c_0\|B\|)^{-M}.\nonumber
\end{align}
Then with constants $c_\nu=\nu^\nu (Mc_0)^{-\nu}$ and $c_{**}=(3+2(Mc_0\|B\|)^2)nMc_Bc_0^2\|B\|$, there holds
\begin{align*}
\E[\|e_{KM}^\delta\|^2]
\leq&\big(2+2^{2\nu}\|B\|c_{**}c_*\big)c^2_\nu K^{-2\nu}\|w\|^2+(2Mc_0+c_{**}c_*)K\bar{\delta}^2.
\end{align*}
\end{theorem}

\begin{remark}\label{rem:con_optimal}
Let $c=c_0\|B\|M$, which implies  $c_{B}=(1-cM^{-1})^{-M}$ and $c_{B,M}=\sum_{i=1}^{M-1}
(1-(1-cM^{-1})^i)^2$, the condition \eqref{cond:sat_svrg} is satisfied whenever
\begin{equation*}
nM^{-2}\leq (1-c_*^{-1})(4+2c^2)^{-1}c_B^{-1}c_{B,M}^{-1},
\end{equation*}
which holds for $M=\mathcal{O}(n^{\frac12})$ and sufficiently small $c=\mathcal{O}(1)$.
It is instructive to compare the conditions ensuring an optimal convergence rate
of SVRG and SGD: SGD requires the condition $c_0=\mathcal{O}(n^{-1})$ \cite{JinZhouZou:2021},
whereas SVRG requires only $M=\mathcal{O}(n^{\frac12})$ and $c=c_0\|B\|M=\mathcal{O}(1)$.
The latter implies $c_0=\mathcal{O}(n^{-\frac12})$ for SVRG.
Since $\mathcal{O}(n^{-\frac12})$ is much larger than $\mathcal{O}(n^{-1})$ when
the data size $n$ is large, {SVRG should perform better for truly large-scale problems.}
\end{remark}

It is known that SGD with an inadvertent choice of the step size schedule can lead to the
undesirable saturation phenomenon, i.e., the convergence rate does not improve with the regularity
index $\nu$ in Assumption \ref{ass:stepsize}(ii), whenever $\nu$ exceeds the critical value $1/2$
\cite{JahnJin:2020,JinZhouZou:2021}. This is attributed to the inherent variance of the stochastic
gradient estimate used by SGD, and one important issue is to overcome the saturation
phenomenon. The next result sheds further insight into this phenomenon by comparing the mean squared error
of the (epochwise) SVRG iterate with that of the corresponding SGD iterate:
it gives a refined comparison between the variance components of SVRG and SGD iterates, in view of the
bias-variance decomposition. In particular,
it shows that the built-in variance reduction mechanism of SVRG does reduce the variance
component of the error, which represents a distinct feature of SVRG over SGD, especially
alleviating the step size restriction for achieving the optimal convergence.

\begin{theorem}\label{thm:svrg-sgd}
Let Assumption \ref{ass:stepsize}(i) and (iii) be fulfilled and the constants $c_0$, $n$ and $M$
satisfy, with the constant $c'_B=(1-c_0\|B\|)^{-2(M-1)}$,
\begin{equation}\label{cond:compare}
(M-1)^2c_0^2\|B\|^2\leq  (2c'_B)^{-1} \mbox{ and }
(M+1)^2\leq  (2c'_B)^{-1}(n-1).
\end{equation}
For any $K\geq 0$, let $R_1$ and $R_2$ be measurable with respect to $\mathcal{F}^c_{KM}$
and $R_1$ is combination of $M_0$ and $H_k$ {\rm(}cf. \eqref{eqn:Hk} for the definition{\rm)}.
Then for $\zeta$ defined in Section \ref{ssec:prelim}, there holds
\begin{align*}
\E[\|R_1 (e_{KM}^\delta-B^{-1}\zeta)+R_2\|^2]\leq\E[\|R_1 (\hat{e}_{KM}^\delta-B^{-1}\zeta)+R_2\|^2].
\end{align*}
\end{theorem}

\begin{remark}\label{cond:M_comp}
Let {$c:=c_0\|B\|(M-1)$}, which implies $c'_B=(1-c(M-1)^{-1})^{-2(M-1)}$. Then
condition \eqref{cond:compare} can be rewritten as
\begin{equation*}
   c^2\leq  2^{-1}(c'_B)^{-1}\quad\mbox{and}\quad (M+1)^2\leq  2^{-1}(c'_B)^{-1}n.
\end{equation*}
The first essentially requires $c<\frac12$. For any $M\geq2$, $c'_B\leq 2 e^{2c}$,
the condition can be satisfied by
$2c e^c \leq 1$ and $M+1\leq 2^{-1}e^{-c}n^{\frac12}$.
\end{remark}

Last we briefly comment on the overall analysis strategy {for proving Theorems \ref{thm:main} and \ref{thm:svrg-sgd}.}
The overall strategy is to derive the recursion of the epochwise SVRG iterate $x_{KM}^\delta$ (and also
the SGD iterate $\hat{x}_{KM}^\delta$), for any  $K=0,1,\cdots$, i.e., at the anchor
points only, and then bound the error $e_{KM}^\delta := x_{KM}^\delta-x^\dag$ by bias-variance
decomposition
\begin{align*}
  \E[\|x_{KM}^\delta - x^\dag\|^2]  & = \|\E[x_{KM}^\delta] - x^\dag\|^2 + \E[\|x_{KM}^\delta-\E[x_{KM}^\delta]\|^2].
\end{align*}
The two terms on the right hand side represent respectively the bias of the error due to early stopping
and data noise and the computational variance of error due to randomness of the gradient estimate.
These are analyzed in Proposition \ref{prop:bias-var} and Lemma \ref{lem:weierr_svrg}, respectively, and allow
proving the convergence rate in Theorem \ref{thm:main}. The analysis of the variance component relies
on a {novel refined decomposition} into terms that are more tractable to estimate for both SVRG and SGD. This
decomposition is also crucial for the comparative study between SVRG and SGD, where a careful componentwise
comparison of the decomposition allows establishing Theorem \ref{thm:svrg-sgd}. Note that the decomposition
relies heavily on the constant step size schedule, and thus the overall analysis differs greatly
from existing analysis of the SGD in the lens of regularization theory \cite{JinLu:2019,JinZhouZou:2020,JinZhouZou:2021} or
the analysis of SGD in statistical learning theory
\cite{YingPontil:2008,TarresYao:2014,DieuleveutBach:2016,LinRosasco:2017,PillaudRudiBach:2018}.
The extension of the analysis to a general step size schedule represents an interesting future research problem.

\section{Error decomposition}\label{sec:decomp}

In this part, we present several preliminary results, especially error decompositions for SVRG and SGD
iterates. The decompositions play a central role in the convergence rates analysis and comparative analysis
in Sections  \ref{sec:conv} and \ref{sec:comp}, respectively.

\subsection{Notation and preliminary estimates}\label{ssec:prelim}

First we introduce several shorthand notation. Below, we denote the SVRG iterates for the exact data $y^\dag$
and noisy data $y^\delta$ by $x_k$ and $x_k^\delta$, respectively, and that for
SGD by $\hat{x}_k$ and $\hat{x}_k^\delta$, respectively. We use extensively the following
shorthand notation for any $k=0,1,\cdots$:
\begin{align*}
&e_k=x_k-x^\dag,\quad e_k^\delta=x_k^\delta-x^\dag,\quad
\hat{e}_k=\hat{x}_k-x^\dag,\quad \hat{e}_k^\delta=\hat{x}_k^\delta-x^\dag,\\
&\bar A= n^{-\tfrac12}A,\quad \bar \xi = n^{-\tfrac12}\xi,\quad \bar{\delta}=n^{-\tfrac12}\delta,\quad M_0=I-c_0B,\quad \zeta={\bar A}^t\bar \xi,\\
&P_k=I-c_0a_{i_k}a_{i_k}^t,\quad N_k=B-a_{i_k}a_{i_k}^t,\quad \zeta_k=a_{i_k} \xi_{i_k}.
\end{align*}
Note that $P_k$ is the random update operator for the iteration, and we have the identity
$P_k=M_0+c_0N_k$ trivially. For all $k\in \mathbb{N}$, let
\begin{equation}\label{eqn:Hk}
H_k= G_{k+1}N_k,
\quad \mbox{with }
  G_k=\left\{\begin{aligned}
    \prod_{i=k}^{k_M+M-1}P_i, &\quad k\neq KM,\\
    I, &\quad k = KM.
  \end{aligned}\right.
\end{equation}
Clearly, $H_{KM-1}=N_{KM-1}.$ By definition, we have the following identity
\begin{align}
  G_{KM+j}&= G_{KM+j+1}P_{KM+j}=G_{KM+j+1}(M_0+c_0N_{KM+j})\nonumber\\
    &=G_{KM+j+1}M_0+c_0H_{KM+j},\quad j=1,\cdots,M-1.\label{eqn:recur-G}
\end{align}
These notations are useful for representing the (epochwise) SVRG iterates $x^\delta_{KM}$, cf. Proposition \ref{prop:bias-var}.
The following simple identity will be used extensively.
\begin{lemma}
The following identity holds
\begin{equation}\label{eqn:G}
  G_{KM+i} = M_0^{M-i} + c_0\sum_{j=0}^{M-i-1}H_{KM+i+j}M_0^j,\quad i =1,\ldots,M-1.
\end{equation}
\begin{proof}
It follows directly from the definition of $G_k$ and $H_k$ and the identity \eqref{eqn:recur-G} that
\begin{align*}
  G_{KM+i} &=G_{KM+i+1}M_0+c_0\sum_{j=0}^0H_{KM+i}M_0^j =  G_{KM+i+2}M_0^2+c_0\sum_{j=0}^1H_{KM+i+j}M_0^j\\
           &=\ldots = M_0^{M-i} + c_0\sum_{j=0}^{M-i-1}H_{KM+i+j}M_0^j.
\end{align*}
This shows the desired identity.
\end{proof}
\end{lemma}

We use extensively the following direct consequence of Assumption \ref{ass:stepsize}(iii).
\begin{lemma}\label{lem:commut}
Under Assumption \ref{ass:stepsize}(iii), the matrices $M_0$, $B$, $P_{i_j}$ and $N_{i_{j'}}$ are commutative for any $j$ and $j'$.
\end{lemma}
\begin{proof}
Note that, for any $j$ and $j'$, we have
\begin{align*}
&B=n^{-1}\sum_{i=1}^n a_i a_i^t, \quad M_0=I-c_0B=I-c_0n^{-1}\sum_{i=1}^n a_i a_i^t,\\
&P_{i_j}=I-c_0a_{i_j}a_{i_j}^t, \quad N_{i_{j'}}=B-a_{i_{j'}}a_{i_{j'}}^t=n^{-1}\sum_{i=1}^n a_i a_i^t-a_{i_{j'}}a_{i_{j'}}^t.
\end{align*}
It suffices to show the claim that $a_i a_i^t$ and $a_j a_j^t$ are commutative for any $i,j=1,\cdots,n$.
This claim is trivial when $i=j$. If $i\neq j$, by Assumption \ref{ass:stepsize}(iii), there holds $a_i^t a_j=0=a_j^t a_i$.
\end{proof}

We also state an identity which is crucial for the proofs of Theorems \ref{thm:main} and \ref{thm:svrg-sgd}.
\begin{lemma}\label{lem:bdd_n}
Let Assumption \ref{ass:stepsize}(iii) be fulfilled. Then for any diagonal matrix $D\in\mathbb{R}^{m\times m}$
and any vector $v\in\mathbb{R}^m$, which are independent of $i_j$, the following identities hold
\begin{align*}
\E[\|VDV^t N_{j}v\|^2]=(n-1)\E[\|VDV^tBv\|^2],\\
\E[\| VDV^t(\zeta_{j}-\zeta)\|^2]=(n-1)\E[\|VDV^t\zeta\|^2].
\end{align*}
\end{lemma}
\begin{proof}
Recall the standard bias-variance decomposition: for any matrix $R$ and filtration $\mathcal{F}_a$,
$$\E[\|R-\E[R|\mathcal{F}_a]\|^2|\mathcal{F}_a]=\E[\|R\|^2|\mathcal{F}_a]-\|\E[R|\mathcal{F}_a]\|^2.$$
Then the identity $N_{j}=B-a_{i_{j}}a_{i_{j}}^t={\E_{j}}[a_{i_{j}}a_{i_{j}}^t]-a_{i_{j}}a_{i_{j}}^t$ gives
\begin{align*}
\E_{j}[\|VDV^t N_{j}v\|^2]=&\E_{j}[\|VDV^ta_{i_{j}}a_{i_{j}}^tv\|^2]-\|  VDV^tBv\|^2\\
=&n^{-1}\sum_{i=1}^n\|   VDV^ta_{i}a_{i}^tv\|^2-\|  VDV^tBv\|^2,
\end{align*}
where $a_ia_i^t v=A^t(a_i^t v)b_i$ with $b_i=(0,\ldots,0,1,0,\ldots,0)^t\in\mathbb{R}^n$ being the ${i}$th canonical Cartesian basis vector.
By Assumption \ref{ass:stepsize}(iii), $DV^tA^t=D\Sigma $ is diagonal, and hence
\begin{align*}
&n^{-1}\sum_{i=1}^n\|   VDV^ta_{i}a_{i}^tv\|^2=n^{-1}\sum_{i=1}^n\| VDV^tA^t (a_{i}^tv)b_{i}\|^2\\
=&n^{-1}\| VDV^tA^t \sum_{i=1}^n (a_{i}^tv)b_{i}\|^2=n\|  VDV^tBv\|^2.
\end{align*}
This shows the first identity. Similarly, since $\E_{j}[\zeta_{j}]=\zeta$, by rewriting
$\zeta_{j}$ as  $\zeta_{j}=a_{i_j} \xi_{i_j}=A^t \xi_{i_j} b_{i_j}$, we obtain the
second identity. This completes the proof of the lemma.
\end{proof}

Next we recall two technical estimates; see the appendix for the proof.
\begin{lemma}\label{lem:kernel}
Let Assumption \ref{ass:stepsize}(i) be fulfilled. For any $s\geq 0$, $t\in[0,1]$ and $K\in \mathbb{N}$,
there hold
\begin{align*}
\|B^{-t} (I-M_0^{KM})\|\leq& (Mc_0)^{t}K^{t}\quad\mbox{and}\quad
\|B^s M_0^{KM}\|\leq s^s (Mc_0)^{-s} K^{-s}:=c_sK^{-s}.
\end{align*}
\end{lemma}

\subsection{Error decomposition}
Now we derive error decompositions for the (epochwise) SVRG error $e_{KM}^\delta
\equiv x_{KM}^\delta-x^\dag$ and the SGD error $\hat e_{KM}^\delta\equiv \hat
x_{KM}^\delta-x^\dag$ into the bias and variance components. These representations follow from direct but
lengthy computation using the definitions the SVRG and SGD iterates, and the
detailed proof is deferred to the appendix.
\begin{proposition}\label{prop:bias-var}
Under Assumption \ref{ass:stepsize}(i), for any $K\geq 1$, there hold
\begin{align*}
\E[e_{(K+1)M}^\delta]=&M_0^{(K+1)M}e_0^\delta+(I-M_0^{(K+1)M})B^{-1} \zeta\\
e_{(K+1)M}^\delta-\E[e_{(K+1)M}^\delta]=&\sum_{j=0}^K M_0^{(K-j)M} L_j(\zeta-B e_{jM}^\delta),
\end{align*}
with the random matrices $L_j$ defined by
\begin{equation*}
  L_j=c_0\sum_{i=1}^{M-1}H_{jM+i}(I-M_0^i)B^{-1}.
\end{equation*}
\end{proposition}

The next result gives an analogous bias-variance decomposition for the SGD iterate $\hat x_{KM}^\delta$.
Note that when compared with Proposition \ref{prop:bias-var}, the expressions of $\E[x_{KM}^\delta]$ and
$\E[\hat x_{KM}^\delta]$ are actually identical,
since both methods use an unbiased estimate for the gradient. Their difference lies in the
variance component, which will be the main focus of the analysis below.
\begin{proposition}\label{prop:bias-var_sgd}
Under Assumption \ref{ass:stepsize}(i), for any $K\geq 0$, $\hat{e}_{(K+1)M}^\delta$ satisfies
\begin{align*}
\E[\hat{e}_{(K+1)M}^\delta]=&M_0^{(K+1)M}\hat{e}_{0}^\delta+(I-M_0^{(K+1)M})B^{-1} \zeta,\\
\hat{e}_{(K+1)M}^\delta-\E[\hat{e}_{(K+1)M}^\delta]
=&c_0^2\sum_{j=0}^K\sum_{i=0}^{M-2}\sum_{t=0}^{M-i-2} M_0^{(K-j)M}H_{jM+i+t+1}M_0^t (\zeta_{jM+i}-\zeta)\\
&+c_0\sum_{j=0}^K\sum_{i=0}^{M-1} M_0^{(K-j)M}\big( H_{jM+i}\big(M_0^i\hat{e}_{jM}^\delta
+(I-M_0^i) B^{-1}\zeta\big)\\
&+ M_0^{M-i-1}(\zeta_{jM+i}-\zeta)\big).
\end{align*}
\end{proposition}

\begin{remark}
Equation \eqref{eqn:svrg-it-KM} in the proof (in the appendix) indicates that at the snapshot point $x_{KM}^\delta$, SVRG performs a
gradient descent step, and in-between the snapshot points, the update direction is a linear combination
between gradient and gradient offset {\rm(}between the current iterate and the anchor point{\rm)}.
Thus in this sense, SVRG is actually a hybridization of the Landweber method and SGD. Note that since $J'(x_{k_M}^\delta)$
is independent of the random index $i_k$ and the gap $f_{i_k}'(x_k^\delta )-f_{i_k}'(x_{k_M}^\delta)$ is
independent of the noise $\xi_{i_k}$ for linear inverse problems, the SVRG iterate $x_k^\delta$ does not
actually depend on $\xi_{i_k}$. This property contributes to the
variance reduction, and constitutes one major difference between SVRG and SGD in terms of the noise influence.
\end{remark}

\section{Proof of Theorem \ref{thm:main}}\label{sec:conv}

Now we prove the convergence rate for SVRG in Theorem \ref{thm:main}.
We begin with bounding the mean squared residual $\E[\|R_1 (e_{KM}^\delta-B^{-1}\zeta)+R_2\|^2]$ and weighted
variance $\E[\|R_1(e_{KM}^\delta-\E[e_{KM}^\delta])\|^2]$, where the quantities $R_1$ and $R_2$ are measurable with
respect to the filtration $\mathcal{F}^c_{KM}$ and commutative with $B$, $M_0$, $\{P_k\}$ and $\{N_k\}$ for any $k\geq 0$.
The specific forms of $R_1$ and $R_2$ arise from the refined decompositions of SVRG errors in Lemma \ref{lem:weierr_svrg}
and SGD errors in Lemma \ref{lem:weierr_sgd}, in order to carry out the componentwise comparison between them;
see the proof of Theorem \ref{thm:svrg-sgd} in Section \ref{sec:comp} for further details.
\begin{lemma}\label{lem:weierr_svrg}
Under Assumption \ref{ass:stepsize}(i) and (iii), for any $K\geq 0$, let $R_1$ and $R_2$ be measurable with
respect to $\mathcal{F}^c_{(K+1)M}$  and commutative with $B$, $M_0$, $\{P_k\}$ and $\{N_k\}$, for any $k\geq 0$. Then there hold
\begin{align*}
\E[\|R_1 (e_{(K+1)M}^\delta-B^{-1}\zeta)+R_2\|^2]=&{\rm I_{0}}+\sum_{j=0}^K{\rm I}_{1,j},\\
\E[\|R_1(e_{(K+1)M}^\delta-\E[e_{(K+1)M}^\delta])\|^2]=&\sum_{j=0}^K{\rm I}_{1,j},
\end{align*}
with the terms ${\rm I}_0$ and ${\rm I}_{1,j}$ given by
\begin{align}
{\rm I}_{0}=&\E[\|R_1M_0^{(K+1)M}(e_{0}^\delta-B^{-1}\zeta)+R_2\|^2],\label{eqn:decom-svrg-0}\\
{\rm I}_{1,j}=&c_0^2\sum_{i=1}^{M-1}\E[\|R_1 M_0^{(K-j)M} H_{jM+i}(I-M_0^i)( e_{jM}^\delta-B^{-1}\zeta)\|^2].\label{eqn:decom-svrg-1j}
\end{align}
\end{lemma}

Now we bound the mean squared (generalized) residual $\E[\|R_1(e_{KM}^\delta-B^{-1}\zeta)\|^2]$ of the
epochwise SVRG iterate $x_{KM}^\delta$. This bound is useful in the proof of Theorem \ref{thm:main}
below. {The proof relies on mathematical induction, and the decomposition in Lemma \ref{lem:weierr_svrg}.}
\begin{theorem}\label{thm:sat_svrg_res}
Let Assumption \ref{ass:stepsize}(i) and (iii) be fulfilled, $R_1$ be a combination of $M_0$ and $B$, and
$c_*>1$ be chosen such that \eqref{cond:sat_svrg} holds.
Then for any $K\geq0$, there holds
\begin{align*}
\E[\|R_1(e_{KM}^\delta-B^{-1}\zeta)\|^2]\leq c_*\|R_1M_0^{\frac{KM}{2}} (e_{0}^\delta-B^{-1}\zeta)\|^2.
\end{align*}
\end{theorem}
\begin{proof}
We prove the theorem by mathematical induction. The case $K=0$ holds true trivially. Now assume that the assertion
holds up to some $K\geq0$, i.e.,
\begin{align}\label{eqn:svrg-hypo}
\E[\|R_1(e_{jM}^\delta-B^{-1}\zeta)\|^2]\leq c_{*}\|R_1M_0^{\frac{jM}{2}} (e_{0}^\delta-B^{-1}\zeta)\|^2,\quad j=0,1,\ldots,K,
\end{align}
and we prove it for the case $K+1$. Lemma \ref{lem:weierr_svrg} with $R_2=0$ gives
\begin{align*}
\E[\|R_1 (e_{(K+1)M}^\delta-B^{-1}\zeta)\|^2]=&{\rm I_{0}}+\sum_{j=0}^K{\rm I}_{1,j},
\end{align*}
with the terms ${\rm I}_0$ and ${\rm I}_{1,j}$ given by \eqref{eqn:decom-svrg-0} (with $R_2=0$) and \eqref{eqn:decom-svrg-1j}.
Note that $V^tR_1 M_0^{(K-j)M}V$ is diagonal, then direct computation with Lemmas \ref{lem:commut} and
\ref{lem:bdd_n}, the inequalities $\|G_{jM+i+1}\|\leq1$ and $\|I-M_0^i\|=1-(1-c_0\|B\|)^i$ and
the definition of the constant $c_{B,M}$ in Theorem \ref{thm:main} gives
\begin{align*}
{\rm I_{0}}\leq&\E[\|R_1M_0^{\frac{(K+1)M}{2}}(e_{0}^\delta-B^{-1}\zeta)\|^2],\\
{\rm I}_{1,j}
\leq&c_0^2\sum_{i=1}^{M-1}\|I-M_0^i\|^2\|G_{jM+i+1}\|^2\E[\|R_1 M_0^{(K-j)M} N_{jM+i}( e_{jM}^\delta-B^{-1}\zeta)\|^2]\\
\leq&nc_0^2c_{B,M}\E[\|R_1 M_0^{(K-j)M} B( e_{jM}^\delta-B^{-1}\zeta)\|^2]\\ 
\leq&nc_0^2c_{B,M}\|M_0^{-\frac{M}{2}}\|^2\|M_0^{\frac{(K-j)M}{2}} B\|^2\E[\|R_1 M_0^{\frac{(K-j+1)M}{2}} ( e_{jM}^\delta-B^{-1}\zeta)\|^2].
\end{align*}
This, the induction hypothesis \eqref{eqn:svrg-hypo}, and the identity
\begin{equation}\label{eqn:cB}
  \|M_0^{-\frac{M}{2}}\|^2=(1-c_0\|B\|)^{-M}:= c_B
\end{equation}
give
\begin{align*}
\sum_{j=0}^{K}{\rm I}_{1,j}\leq nc_0^2c_Bc_{B,M}c_*\sum_{j=0}^K\|M_0^{\frac{(K-j)M}{2}}B\|^2\|R_1M_0^{\frac{(K+1)M}{2}}(e_{0}^\delta-B^{-1}\zeta)\|^2.
\end{align*}
By Lemma \ref{lem:kernel},
$$\|M_0^{\frac{(K-j)M}{2}}B\|\leq 2((K-j)Mc_0)^{-1},\quad j=0,\cdots,K-2,$$
and consequently,
\begin{equation}\label{eqn:sum-M0B}
  \sum_{j=0}^K\|M_0^{\frac{(K-j)M}{2}}B\|^2\leq2\|B\|^2+ 4c_0^{-2}M^{-2}\sum_{j=0}^{K-2}(K-j)^{-2} \leq (4+2(Mc_0\|B\|)^2)c_0^{-2}M^{-2}.
\end{equation}
The preceding estimates together imply
\begin{align*}
   &\E[\|R_1(e_{(K+1)M}^\delta-B^{-1}\zeta)\|^2]\\
 \leq& \big(1+(4+2(Mc_0\|B\|)^2)nM^{-2}c_Bc_{B,M}c_*\big)\|R_1M_0^{\frac{(K+1)M}{2}}(e_{0}^\delta-B^{-1}\zeta)\|^2.\nonumber
\end{align*}
The condition on $c_*$ from \eqref{cond:sat_svrg} shows the induction step, and this completes
the proof of the theorem.
\end{proof}

Setting $R_1=n^{\tfrac12} B^{\tfrac12}$ in Theorem \ref{thm:sat_svrg_res} gives
an upper bound on the mean squared residual $\E[\|Ax_{KM}^\delta-y^{\delta}\|^2]$ of the
(epochwise) SVRG iterate $x_{KM}^\delta$. {Note that the mean squared residual
consists of one decaying term related to the source condition in Assumption \ref{ass:stepsize}(ii)
and one constant term related to the noise level. In particular, it is essentially bounded,
independent of the iteration index. This behavior is similar to that for the standard Landweber method.}
\begin{corollary}\label{cor:sat_svrg_res}
Under Assumption \ref{ass:stepsize} and condition \eqref{cond:sat_svrg}, there holds
\begin{align*}
\E[\|Ax_{KM}^\delta-y^\delta\|^2]\leq2^{2\nu+2}c^2_{\nu+\frac12}nc_{*} K^{-2\nu-1}\|w\|^2+2nc_{*}\bar\delta^2.
\end{align*}
\end{corollary}
\begin{proof}
Theorem \ref{thm:sat_svrg_res} and the triangle inequality imply (noting $e_0^\delta=e_0$)
\begin{align*}
\E[\|Ax_{KM}^\delta-y^\delta\|^2]=&\E[\|n^{\frac12}B^{\frac12}(e_{KM}^\delta-B^{-1}\zeta)\|^2]
\leq nc_*\|B^{\frac12}M_0^{\frac{KM}{2}} (e_{0}^\delta-B^{-1}\zeta)\|^2,\\
\leq &  2nc_*\|M_0^{\frac{KM}{2}}B^{\frac12} e_{0}^\delta\|^2+2nc_*\|M_0^{\frac{KM}{2}}B^{-\frac12}\zeta\|^2.
\end{align*}
Meanwhile, it follows from Lemma \ref{lem:kernel} and the source condition in Assumption \ref{ass:stepsize}(ii) that
\begin{align*}
\|M_0^{\frac{KM}{2}}B^{\frac12} e_{0}\|\leq & 2^{\nu+\frac12}c_{\nu+\frac12} K^{-\nu-\frac12}\|w\|,\\
\|M_0^{\frac{KM}{2}}B^{-\frac12}\zeta\|^2\leq&\|M_0^{\frac{KM}{2}}B^{-\frac12}\bar{A}^t\|^2\|\bar{\xi}\|^2\leq \bar{\delta}^2.
\end{align*}
Combining the preceding estimates gives the desired assertion.
\end{proof}

Now we can present the proof of Theorem \ref{thm:main}. {The proof employs the representation
in Theorem \ref{thm:sat_svrg_res}, and follows by directly bounding the involved terms using Lemma
\ref{lem:kernel} (under Assumption \ref{ass:stepsize}(ii)) and Theorem \ref{thm:sat_svrg_res}.}
\begin{proof}
By Lemma \ref{lem:weierr_svrg}, setting $R_1=I$ and $R_2=B^{-1}\zeta$ gives
\begin{align*}
\E[\|e_{(K+1)M}^\delta\|^2]\leq{\rm I_{0}}+\sum_{j=0}^K{\rm I}_{1,j},
\end{align*}
with the terms ${\rm I}_0$ and ${\rm I}_{1,j}$ given by \eqref{eqn:decom-svrg-0}
and \eqref{eqn:decom-svrg-1j}, respectively. Now we bound them separately. By the triangle inequality,
Assumption \ref{ass:stepsize}(ii) and Lemma \ref{lem:kernel}, we deduce
\begin{align*}
{\rm I_{0}}&=\|M_0^{(K+1)M}e_{0}+(I-M_0^{(K+1)M})B^{-1}\zeta\|^2\\
&\leq 2 \|M_0^{(K+1)M}e_{0}\|^2+2\|(I-M_0^{(K+1)M})B^{-1}\bar{A}^t\bar{\xi}\|^2\\
&\leq 2c^2_\nu (K+1)^{-2\nu}\|w\|^2+2Mc_0(K+1)\bar{\delta}^2.
\end{align*}
Meanwhile, \eqref{eqn:decom-svrg-1j} with $R_1=I$ gives
\begin{align*}
{\rm I}_{1,j}=c_0^2\sum_{i=1}^{M-1}\E[\| M_0^{(K-j)M} (I-M_0^i)H_{jM+i}( e_{jM}^\delta-B^{-1}\zeta)\|^2].
\end{align*}
Note that by Lemma \ref{lem:commut}, the matrices $I-M_0^i$ and $H_{jM+i}$ are commuting, and
$H_{jM+i}=G_{jM+i+1}N_{jM+i}$. Thus by Lemma \ref{lem:bdd_n} (with $V^tM_0^{(K-j)M} (I-M_0^i)G_{jM+i+1}V$
being diagonal) and $\|G_{jM+i+1}\|\leq1$, we obtain
\begin{align*}
{\rm I}_{1,j}= &(n-1)c_0^2\sum_{i=1}^{M-1}\|  M_0^{(K-j)M} (I-M_0^i) G_{jM+i+1}B( e_{jM}^\delta-B^{-1}\zeta)\|^2\\
\leq& nc_0^2\sum_{i=1}^{M-1}\E[\| M_0^{(K-j)M} B(I-M_0^i)( e_{jM}^\delta-B^{-1}\zeta)\|^2].
\end{align*}
Next by the identity
$$c_0\sum_{i=0}^{j-1}M_0^i=(I-M_0^{j})B^{-1},$$
 the trivial inequality
$(\sum_{t=0}^{i-1}a_t)^2\leq i\sum_{t=0}^{i-1}a_t^2$, and $\|M_0\|\leq 1$, we have
\begin{align*}
{\rm I}_{1,j}\leq&nc_0^4\sum_{i=1}^{M-1}\E[\| M_0^{(K-j)M} B^2\sum_{t=0}^{i-1}M_0^t( e_{jM}^\delta-B^{-1}\zeta)\|^2]\\
\leq&nc_0^4\sum_{i=1}^{M-1}i\sum_{t=0}^{i-1}\E[\| M_0^{(K-j)M} BM_0^t( Be_{jM}^\delta-\zeta)\|^2]\\
\leq&nc_0^4\sum_{i=1}^{M-1}i^2\E[\| M_0^{(K-j)M} B( Be_{jM}^\delta-\zeta)\|^2].
\end{align*}
Since $\sum_{i=1}^{M-1}i^2\leq3^{-1} M^3$, it follows from Theorem \ref{thm:sat_svrg_res} and \eqref{eqn:cB} that 
\begin{align*}
{\rm I}_{1,j}\leq&3^{-1}nM^3c_0^4\E[\| M_0^{(K-j)M} B( Be_{jM}^\delta-\zeta)\|^2].\\
\leq&3^{-1}nM^3c_0^4 \|M_0^{-\frac{M}{2}}\|^2\|M_0^{\frac{(K-j)M}{2}}B\|^2\E[\|M_0^{\frac{(K-j+1)M}{2}}(Be_{jM}^\delta-\zeta)\|^2]\\
\leq& 3^{-1}nc_BM^3c_0^4c_*\|M_0^{\frac{(K-j)M}{2}}B\|^2\|M_0^{\frac{(K+1)M}{2}}(Be_{0}^\delta-\zeta)\|^2.
\end{align*}
This and the inequality \eqref{eqn:sum-M0B} imply
\begin{align*}
  \sum_{j=0}^{K}{\rm I}_{1,j} &\leq 3^{-1}(4+2(Mc_0\|B\|)^2)nMc_Bc_0^2c_*\|M_0^{\frac{(K+1)M}{2}}(Be_{0}^\delta-\zeta)\|^2\\
    &\leq (3+2(Mc_0\|B\|)^2)nMc_Bc_0^2c_*\big(2^{2\nu}\|B\|^2c^2_\nu(K+1)^{-2\nu}\|w\|^2+\|B\|\bar{\delta}^2\big).
\end{align*}
The last two estimates together yield
\begin{align*}
\E[\|e_{(K+1)M}^\delta\|^2]
\leq&\big(2+2^{2\nu}(3+2(Mc_0\|B\|)^2)nMc_Bc_0^2\|B\|^2c_*\big)c^2_\nu (K+1)^{-2\nu}\|w\|^2
\\
&+(2Mc_0+(3+2(Mc_0\|B\|)^2)nMc_Bc_0^2\|B\|c_*)(K+1)\bar{\delta}^2.
\end{align*}
This completes the proof of the theorem.
\end{proof}

\section{Proof of Theorem \ref{thm:svrg-sgd}}\label{sec:comp}
This section is devoted to the proof of Theorem \ref{thm:svrg-sgd}, and presents a comparative
study on the variance $\E[\|e_{KM}^\delta-\E[e_{KM}^\delta]\|^2]$ of SVRG iterates with $\E[\|\hat
e_{KM}^\delta-\E[\hat e_{KM}^\delta]\|^2]$ of SGD iterates. First we give a bias-variance
decomposition of the SGD iterate $\hat x_{KM}^\delta$, in analogy with Lemma \ref{lem:weierr_svrg}.
{The representations in Lemmas \ref{lem:weierr_svrg} and \ref{lem:weierr_sgd} facilitate
the comparison between the variance components directly, which, under certain conditions, enables
comparing the variance of SVRG and SGD iterates.}
\begin{lemma}\label{lem:weierr_sgd}
Under Assumption \ref{ass:stepsize}(i) and (iii), for any $K\geq 0$, let $R_1$ and $R_2$ be measurable with
respect to $\mathcal{F}^c_{(K+1)M}$  and commutative with $B$, $M_0$, $\{P_k\}$ and $\{N_k\}$, for any $k\geq 0$. Then there hold
\begin{align*}
\E[\|R_1 (\hat{e}_{(K+1)M}^\delta-B^{-1}\zeta)+R_2\|^2]=&{\rm I_0}+\sum_{j=0}^K({\rm I}_{2,j}+{\rm I}_{3,j}),\\
\E[\|R_1(\hat{e}_{(K+1)M}^\delta-\E[\hat{e}_{(K+1)M}^\delta])\|^2]=&\sum_{j=0}^K({\rm I}_{2,j}+{\rm I}_{3,j}),
\end{align*}
with ${\rm I}_0$ given by \eqref{eqn:decom-svrg-0} and ${\rm I}_{2,j}$ and ${\rm I}_{3,j}$ given by
\begin{align}
{\rm I}_{2,j}
=&c_0^2\sum_{i=0}^{M-1}\E[\|R_1 M_0^{(K-j)M}\big( H_{jM+i}M_0^i(\hat{e}_{jM}^\delta-B^{-1}\zeta)
+H_{jM+i}B^{-1}\zeta+ M_0^{M-i-1}(\zeta_{jM+i}-\zeta)\big)\|^2],\label{eqn:decom-sgd-2j}\\
{\rm I}_{3,j}=&c_0^4 \sum_{i=1}^{M-1}\sum_{t=0}^{i-1}\E[\|R_1 M_0^{(K-j)M}H_{jM+i}M_0^t (\zeta_{jM+i-1-t}-\zeta)\|^2].\label{eqn:decom-sgd-3j}
\end{align}
\end{lemma}

Now, we can prove Theorem \ref{thm:svrg-sgd}. This result states that the variance component of the
SVRG iterate $x_{KM}^\delta$ is indeed smaller than that of the SGD iterate $\hat x_{KM}^\delta$,
as one may expect from the construction of variance reduction, and thus the variance reduction step
does reduce the variance of the iterate, {thereby alleviating the deleterious effect of
the stochastic iteration noise on the convergence of the SVRG iterates. The proof relies heavily on
the explicit representations of the variances for the iterates $x_{KM}^\delta$ and $\hat x_{KM}^\delta$
derived in Lemmas \ref{lem:weierr_svrg} and \ref{lem:weierr_sgd}, and employs mathematical induction,
certain independence relations (cf. \eqref{eqn:bdd_1}--\eqref{eqn:bdd_3}) as well as lengthy computation.}
\begin{proof}
Recall that the assumption on $R_1$ implies that it is commutative with $B$, $M_0$,
$\{P_k\}$ and $\{N_k\}$ for any $k\geq 0$, and that in the inequality, $R_1$ and $R_2$
are measurable with respect to $\mathcal{F}^c_{jM}$ (when considering $e^\delta_{jM}$). These facts will be
used extensively without explicit mentioning below. The proof proceeds by mathematical
induction. The case $K=0$ is trivial since $\hat{e}_0^\delta=
e_0^\delta$. Now suppose that the assertion holds up to some $K$, i.e.,
\begin{align}\label{eqn:svrg-sgd-comp}
\E[\|R_1 (e_{jM}^\delta-B^{-1}\zeta)+R_2\|^2]\leq
\E[\|R_1 (\hat{e}_{jM}^\delta-B^{-1}\zeta)+R_2\|^2],\quad j=0,1,\ldots, K,
\end{align}
and we prove it for $j=K+1$. By Lemmas \ref{lem:weierr_svrg} and \ref{lem:weierr_sgd}, we deduce
\begin{align*}
\E[\|R_1 (e_{(K+1)M}^\delta-B^{-1}\zeta)+R_2\|^2]=&{\rm I_{0}}+\sum_{j=0}^K{\rm I}_{1,j},\\
\E[\|R_1 (\hat{e}_{(K+1)M}^\delta-B^{-1}\zeta)+R_2\|^2]=&{\rm I_0}+\sum_{j=0}^K({\rm I}_{2,j}+{\rm I}_{3,j}),
\end{align*}
with the terms ${\rm I}_{1,j}$, ${\rm I}_{2,j}$ and ${\rm I}_{3,j}$ are given by
\eqref{eqn:decom-svrg-1j}, \eqref{eqn:decom-sgd-2j} and \eqref{eqn:decom-sgd-3j},
respectively. Thus, it suffices to show
\begin{equation}\label{eqn:svrg-sgd-ind}
  {\rm I}_{1,j}\leq {\rm I}_{2,j}+{\rm I}_{3,j},\quad j=0,1,\cdots,K.
\end{equation}
By the inequality $(\sum_{t=1}^ia_i)^2\leq i\sum_{t=1}^ia_i^2$, \eqref{eqn:sum-M}
and the identity $\|M_0^{-1}\|=(1-c_0\|B\|)^{-1}$, we have
\begin{align*}
{\rm I}_{1,j}=&c_0^4\sum_{i=1}^{M-1}\E[\|R_1 M_0^{(K-j)M} H_{jM+i}B\sum_{t=0}^{i-1}M_0^t( e_{jM}^\delta-B^{-1}\zeta)\|^2]\\
\leq&c_0^4\sum_{i=1}^{M-1}i\sum_{t=0}^{i-1}\|M_0^{-i}\|^2\E[\|M_0^i R_1 M_0^{(K-j)M} H_{jM+i}BM_0^t( e_{jM}^\delta-B^{-1}\zeta)\|^2]\\
\leq&c_0^4\sum_{i=1}^{M-1}i(1-c_0\|B\|)^{-2i}\sum_{t=0}^{i-1}\E[\|M_0^i R_1 M_0^{(K-j)M} H_{jM+i}M_0^tB( e_{jM}^\delta-B^{-1}\zeta)\|^2]\\
\leq&c_0^4\sum_{i=1}^{M-1}i(1-c_0\|B\|)^{-2i}\sum_{t=0}^{i-1}\E[\|M_0^i R_1 M_0^{(K-j)M} H_{jM+i}M_0^tB(\hat e_{jM}^\delta-B^{-1}\zeta)\|^2],
\end{align*}
where the last step is due to the induction hypothesis \eqref{eqn:svrg-sgd-comp}. Then
by Lemma \ref{lem:commut}, adding and subtracting suitable terms, and the triangle
inequality, since $\|M_0\|\leq 1$, we deduce (with shorthand notation $c'_B=(1-c_0\|B\|)^{-2(M-1)}$)
\begin{align*}
{\rm I}_{1,j}
\leq&c_0^4\sum_{i=1}^{M-1}i(1-c_0\|B\|)^{-2i}\sum_{t=0}^{i-1}\E[\| R_1 M_0^{(K-j)M+t}B\big( H_{jM+i}M_0^i(\hat e_{jM}^\delta-B^{-1}\zeta)+H_{jM+i}B^{-1}\zeta\\
&\qquad+ M_0^{M-i-1}(\zeta_{jM+i}-\zeta)\big)-R_1 M_0^{(K-j)M+t}\big( H_{jM+i}\zeta+ M_0^{M-i-1}B(\zeta_{jM+i}-\zeta)\big)\|^2]\\
\leq&2(M-1)^2\|B\|^2c'_Bc_0^2{\rm I}_{2,j}+c_0^4\sum_{i=1}^{M-1}i(1-c_0\|B\|)^{-2i}\sum_{t=0}^{i-1}\Big(4\E[\| R_1 M_0^{(K-j)M+t} H_{jM+i}\zeta\|^2]\\
&\qquad+4\E[\| R_1 M_0^{(K-j+1)M+t-i-1} B(\zeta_{jM+i}-\zeta)\|^2]\Big).
\end{align*}
Now Assumption \ref{ass:stepsize}(iii) and the condition on $R_1$ imply that
$V^tR_1 M_0^{s_1}G_{k+1} N_k^{s_3}B^{s_2}V$ is diagonal for any $s_1,s_2\geq0$, $s_3=0,1$ and $k\in\mathbb{N}$. Thus, by Lemma \ref{lem:bdd_n}, we obtain
\begin{align}\label{eqn:bdd_1}
\E[\| R_1 M_0^{(K-j)M+t} H_{jM+i}\zeta\|^2]=&(n-1)\E[\| R_1 M_0^{(K-j)M+t}G_{jM+i+1} B\zeta\|^2],\\
\E[\| R_1 M_0^{(K-j+1)M+t-i-1} B(\zeta_{jM+i}-\zeta)\|^2]=&(n-1)\E[\| R_1 M_0^{(K-j+1)M+t-i-1} B\zeta\|^2],\label{eqn:bdd_2}\\
\E[\|R_1 M_0^{(K-j)M+t}H_{jM+i} (\zeta_{jM+i-1-t}-\zeta)\|^2]=&(n-1)\E[\|R_1 M_0^{(K-j)M+t}H_{jM+i} \zeta\|^2]\nonumber\\
=&(n-1)^2\E[\|R_1 M_0^{(K-j)M+t} G_{jM+i+1}B\zeta\|^2].\label{eqn:bdd_3}
\end{align}
Using the relation $H_{jM+M-1}=N_{jM+M-1}$ and \eqref{eqn:bdd_3} leads to
\begin{align*}
{\rm I}_{3,j} 
=&c_0^4 \sum_{i=1}^{M-2}\sum_{t=0}^{i-1}\E[\|R_1 M_0^{(K-j)M+t}H_{jM+i} (\zeta_{jM+i-1-t}-\zeta)\|^2]\\
&+c_0^4 \sum_{t=0}^{M-2}\E[\|R_1 M_0^{(K-j)M+t}N_{jM+M-1} (\zeta_{jM+M-2-t}-\zeta)\|^2]\\
=&(n-1)c_0^4 \sum_{i=1}^{M-2}\sum_{t=0}^{i-1}\E[\|R_1 M_0^{(K-j)M+t}H_{jM+i} \zeta\|^2]\\
 &+(n-1)^2c_0^4 \sum_{t=0}^{M-2}\E[\|R_1 M_0^{(K-j)M+t}B \zeta\|^2].
\end{align*}
Let ${\rm II}_{j,i,t}=\E[\|R_1 M_0^{(K-j)M+t}H_{jM+i} \zeta\|^2]$, and ${\rm II}_{j,0,t}=\E[\|R_1 M_0^{(K-j)M+t}B \zeta\|^2]$.
Similarly, with the identities \eqref{eqn:bdd_1} and \eqref{eqn:bdd_2}, we deduce
\begin{align*}
{\rm I}_{1,j}
\leq&2(M-1)^2\|B\|^2{c'_B}c_0^2{\rm I}_{2,j}+4c_0^4\sum_{i=1}^{M-2}i(1-c_0\|B\|)^{-2i}\sum_{t=0}^{i-1}\E[\| R_1 M_0^{(K-j)M+t} H_{jM+i}\zeta\|^2]\\
&+4{c'_B}c_0^4(M-1)\sum_{t=0}^{M-2}\E[\| R_1 M_0^{(K-j)M+t} N_{jM+M-1}\zeta\|^2]\\
&+4c_0^4\sum_{i=1}^{M-1}i(1-c_0\|B\|)^{-2i}\sum_{t=0}^{i-1}\E[\| R_1 M_0^{(K-j+1)M+t-i-1} B(\zeta_{jM+i}-\zeta)\|^2]\\
\leq&2(M-1)^2\|B\|^2{c'_B}c_0^2{\rm I}_{2,j}+4(M-2)c_Bc_0^4\sum_{i=1}^{M-2}\sum_{t=0}^{i-1}{\rm II}_{j,i,t}\\
&+4(n-1)(M-1){c'_B}c_0^4\sum_{t=0}^{M-2}\E[\| R_1 M_0^{(K-j)M+t} B\zeta\|^2]\\
&+4(n-1){c'_B}c_0^4\sum_{i=1}^{M-1}i\sum_{t=0}^{i-1}\E[\| R_1 M_0^{(K-j+1)M+t-i-1} B\zeta\|^2].
\end{align*}
Note that $\|M_0^{M-i-1}\|^2\leq1$ for any $1\leq i\leq M-1$. The last two terms on
the right hand side of the inequality, denoted by ${\rm II}$, can be bounded by
\begin{align*}
{\rm II}\leq& 4(n-1){c'_B}c_0^4\Big((M-1)\sum_{t=0}^{M-2}+\sum_{i=1}^{M-1}i\sum_{t=0}^{i-1}\Big)\E[\| R_1 M_0^{(K-j)M+t} B\zeta\|^2]\\
=&4(n-1){c'_B}c_0^4\sum_{t=0}^{M-2}\Big(M-1+\sum_{i=t+1}^{M-1}i\Big){\rm II}_{j,0,t}\\
\leq &2(n-1)(M+1)^2{c'_B}c_0^4 \sum_{t=0}^{M-2}{\rm II}_{j,0,t},
\end{align*}
since $M-1+\sum_{i=t+1}^{M-1}i \leq\frac12 (M+1)^2$, for $ 0\leq t\leq M-2.$ Consequently,
\begin{align*}
{\rm I}_{1,j}\leq&2(M-1)^2\|B\|^2c'_Bc_0^2{\rm I}_{2,j}+4(M-2)c'_Bc_0^4\sum_{i=1}^{M-2}\sum_{t=0}^{i-1}{\rm II}_{j,i,t}\\
   &+2(n-1)(M+1)^2c'_Bc_0^4 \sum_{t=0}^{M-2}{\rm II}_{j,0,t}.
\end{align*}
Now the condition \eqref{cond:compare} implies \eqref{eqn:svrg-sgd-ind}, which shows the induction
step and completes the proof of the theorem.
\end{proof}

\begin{remark}\label{rem:sat_0}
For exact data, i.e., $\delta=0$, $\zeta=0$, $\zeta_i=0$ for any $i\geq0$, the comparative
analysis can be greatly simplified. Indeed, setting $R_1=I$ and $R_2=0$ in the analysis leads to
\begin{align*}
\E[\|e_{(K+1)M}\|^2]\leq{\rm I_{0}}+\sum_{j=0}^K{\rm I}_{1,j},
\end{align*}
with
\begin{align*}
{\rm I_{0}}=\|M_0^{(K+1)M}e_{0}\|^2\quad\mbox{and}\quad
{\rm I}_{1,j}=c_0^2\sum_{i=1}^{M-1}\E[\| M_0^{(K-j)M} H_{jM+i}(I-M_0^i) e_{jM}\|^2].
\end{align*}
Straightforward computation with Lemma \ref{lem:bdd_n} gives
\begin{align*}
{\rm I}_{1,j}\leq&(n-1)c_0^4\sum_{i=1}^{M-1}i^2\E[\| M_0^{(K-j)M} G_{jM+i+1} B^2 e_{jM}\|^2]\\
\leq&(n-1)(M-1)^2c_B'c_0^4\|B\|^2\sum_{i=1}^{M-1}\E[\| M_0^{(K-j)M+i} G_{jM+i+1} B e_{jM}\|^2].
\end{align*}
Similarly, Lemma \ref{lem:weierr_sgd} with $R_1=I$ and $R_2=0$ implies
\begin{align*}
\E[\|\hat{e}_{(K+1)M}\|^2]=&{\rm I_0}+\sum_{j=0}^K{\rm I}_{2,j},
\end{align*}
with
\begin{align*}
{\rm I}_{2,j}&=c_0^2\sum_{i=0}^{M-1}\E[\| M_0^{(K-j)M} H_{jM+i}M_0^i\hat{e}_{jM}\|^2]\\
 &=(n-1)c_0^2\sum_{i=0}^{M-1}\E[\| M_0^{(K-j)M+i} G_{jM+i+1}B\hat{e}_{jM}
\|^2].
\end{align*}
When $c_0\|B\|(M-1)\leq (1-c_0\|B\|)^{(M-1)},$ the conditions for the optimal convergence rate of SVRG
is weaker than that of SGD. With $c=c_0\|B\|(M-1)$ and $c_1=(1-c(M-1)^{-1})^{(M-1)}$,
the conditions can be satisfied if $c\leq c_{1}.$ This short analysis clearly shows the beneficial effect
of variance reduction on the variance of the iterates $x_k^\delta$, and hence SVRG allows larger step
size while maintaining the optimal convergence.
\end{remark}

\section{Numerical experiments and discussions}\label{sec:numer}
In this section, we provide numerical experiments to complement the theoretical findings in Section \ref{sec:main}.
{The experimental setting is identical with that in \cite{JinZhouZou:2021}}. Specifically, we employ three
academic examples, i.e., \texttt{s-phillips} (mildly ill-posed), \texttt{s-gravity} (severely ill-posed) and \texttt{s-shaw}
(severely ill-posed), generated from \texttt{phillips}, \texttt{gravity} and \texttt{shaw}, taken from the \texttt{MATLAB}
package Regutools \cite{P.C.Hansen2007} (available at \url{http://people.compute.dtu.dk/pcha/Regutools/}, last accessed
on August 20, 2020), all of size $n = m = 1000$. To explicitly control the regularity index $\nu$ in Assumption
\ref{ass:stepsize}(ii), we generate $x^\dag$ by $x^\dag = \|(A^tA)^\nu x_e\|_{\ell^\infty}^{-1}(A^tA)^\nu x_e$,
where $x_e$ is the exact solution given by the package, and $\|\cdot\|_{\ell^\infty}$ denotes the Euclidean maximum norm.
The index $\nu$ in Assumption \ref{ass:stepsize}(ii) is slightly larger than the one defined above. The corresponding exact data $y^\dag$ is given by $y^\dag=A x^\dag$ and the noise data $y^\delta$ generated by
\begin{equation*}
  y^\delta_i:=y^\dag_i+\epsilon\|y^\dag\|_{\ell^\infty}\xi_i,\quad i=1,\cdots,n,
\end{equation*}
where $\xi_i$s follow the standard Gaussian distribution, and $\epsilon > 0$ is the relative noise
level. The maximum number of epochs is fixed at $9$e5, where one epoch refers to $\tfrac{nM}{n+M}$ SVRG
iterations or $n$ SGD iterations so that the computational complexity of each method is comparable.
All statistical quantities are computed from 100 runs. We present also numerical results for
the Landweber method (LM) \cite[Chapter 6]{EnglHankeNeubauer:1996} (with a step size $\|A\|^{-2}$),
since it enjoys order optimality. All methods are initialized with $x_0 = 0$.

The accuracy of the reconstructions is measured by the mean squared errors $e_{\rm svrg}=\E[\|x_{k_*}^\delta-x^\dag\|^2]$,
$e_{\rm sgd}=\E[\|\hat x_{k_*}^\delta-x^\dag\|^2]$ for SVRG and SGD, respectively,
and the squared error $e_{\rm lm}=\|x_{k_*}^\delta-x^\dag\|^2$ for LM. The
stopping index $k_*$ (measured in epoch count) is taken such that
the error is smallest along the respective iteration trajectory, due to
a lack of rigorous \textit{a posteriori} stopping rules for SVRG and SGD (the discrepancy principle
is indeed convergent for SGD, without a rate \cite{JahnJin:2020}). The
constant $c$ in the step size $c_0$ is $c=(\max_i(\|a_i\|^2)
)^{-1}$, so that $c_0=\mathcal{O}(cM^{-1})$ for SVRG and $c_0=\mathcal{O}(cn^{-1})$ for SGD.

\subsection{Numerical results for general $A$}

\begin{table}[htp!]
  \centering\small
  \begin{threeparttable}
  \caption{Comparison between SVRG (with $M=100$), SGD and LM for \texttt{s-phillips}.\label{tab:phil}}
    \begin{tabular}{cccccccccccc}
    \toprule
    \multicolumn{2}{c}{Method}&
    \multicolumn{3}{c}{SVRG}&\multicolumn{3}{c}{SGD}&\multicolumn{2}{c}{LM}\\
    \cmidrule(lr){3-5} \cmidrule(lr){6-8} \cmidrule(lr){9-10}
    $\nu$& $\epsilon$ &$c_0$&$e_{\rm svrg}$&$k_{\rm svrg}$&$c_0$&$e_{\rm sgd}$&$k_{\rm sgd}$&$e_{\rm lm}$&$k_{\rm lm}$\\
    \midrule
    $0$&1e-3 & $5c/M$ &  1.67e-2& 4134.35   & $4c/n$  & 1.66e-2 &  4691.28 & 1.65e-2 &  5851
    \cr
       &1e-2 & $5c/M$ & 1.31e-1 &180.95 & $4c/n$  & 1.29e-1 &  204.90  & 1.28e-1 &249\cr
       &5e-2 & $5c/M$ & 5.42e-1  & 96.25    & $4c/n$  & 5.42e-1 &  108.90  & 5.34e-1 &  136    \cr
    \hline
    $1$&1e-3 & $1.5c/M$ & 3.31e-4  & 430.65 & $c/n$   & 3.48e-4 &  539.19  & 2.28e-4 &  157   \cr
    &1e-2 & $1.5c/M$ & 5.96e-3 &41.25    & $c/n$   & 6.64e-3 &57.81     & 5.12e-3 &16
       \cr
       &5e-2 & $1.5c/M$ & 3.22e-2 & 21.45   & $c/n$   & 3.52e-2 &  29.40   & 3.16e-2 &  8
       \cr
    \hline
    $2$&1e-3 & $c/(2M)$& 7.16e-5 & 155.10  &$c/(30n)$& 7.02e-5 &  2115.54 & 3.22e-5 &  19    \cr
       &1e-2 & $c/(2M)$& 1.07e-3 &68.75      &$c/(30n)$& 1.09e-3 &938.70    & 9.82e-4 &8      \cr
       &5e-2 & $c/(2M)$& 2.90e-2& 46.75     &$c/(30n)$& 2.92e-2 &  636.51  & 1.57e-2 &  5    \cr
       \hline
    $4$&1e-3 & $c/(5M)$& 3.05e-5 & 202.95   &$c/(30n)$&  9.77e-5 &  1966.38& 1.30e-5  & 8    \cr
       &1e-2 & $c/(5M)$& 2.41e-3 &142.45    &$c/(30n)$&  2.56e-3 &785.94   & 1.42e-3  &5      \cr
       &5e-2 & $c/(5M)$& 5.20e-2 & 110.00   &$c/(30n)$&  5.23e-2 & 596.73  & 2.49e-2  & 3     \cr
    \bottomrule
    \end{tabular}
    \end{threeparttable}
\end{table}

The numerical results for the three examples with different regularity index $\nu$ and
different noise levels are shown in Tables \ref{tab:phil}--\ref{tab:shaw}, where the
employed constant step size is determined in order to achieve optimal convergence (while maintaining good
computational efficiency). For each fixed regularity index $\nu$, all the errors $e_{\rm svrg}$, $e_{\rm sgd}$ and
$e_{\rm lm}$ decrease to zero as the (relative) noise level $\epsilon$ tends to zero with a certain rate,
and the precise convergence rate depends on the index $\nu$ roughly as the theoretical prediction
$\mathcal{O}(\delta^\frac{4\nu}{2\nu+1})$ (cf. Theorem \ref{thm:main} for SVRG, and
Remark \ref{rem:con_optimal} for SGD). Generally a larger $\nu$ leads to a faster convergence
with respect to $\delta$ as the theory indicates, but the required number of
iterations to reach the optimal error may not necessarily decrease, due to the use of smaller step
sizes. The latter contrasts sharply with that for LM, for which a smoother exact solution $x^\dag$ requires
fewer iterations to reach optimal accuracy (when $\delta$ is fixed). Note that for both SVRG and SGD,
optimal convergence holds only for a sufficiently small step size, and otherwise they
suffer from the undesirable saturation phenomenon, i.e., the error decay may saturate when the index $\nu$
exceeds a certain value, which also concurs with the observation for SGD in
\cite{JahnJin:2020,JinZhouZou:2021}.

\begin{table}[htp!]
  \centering\small
  \begin{threeparttable}
  \caption{Comparison between SVRG (with $M=100$), SGD and LM for \texttt{s-gravity}.\label{tab:gravity}}
    \begin{tabular}{ccccccccccccccc}
    \toprule
    \multicolumn{2}{c}{Method}&
    \multicolumn{3}{c}{SVRG}&\multicolumn{3}{c}{SGD}&\multicolumn{2}{c}{LM}\cr
    \cmidrule(lr){3-5} \cmidrule(lr){6-8}
    \cmidrule(lr){9-10}
    $\nu$& $\epsilon$ &$c_0$&$e_{\rm svrg}$&$k_{\rm svrg}$&$c_0$&$e_{\rm sgd}$&$k_{\rm sgd}$&$e_{\rm lm}$&$k_{\rm lm}$\\
    \midrule
    $0$&1e-3  & $c/10$ & 9.50e-2 &  5495.05 & $c/20$  & 9.37e-2&  1000.50& 9.39e-2  & 27201 \cr
       &1e-2  & $c/10$ & 5.98e-1  &217.80     & $c/20$  & 5.81e-1 &34.11     & 5.73e-1   &793    \cr
       &5e-2  & $c/10$ & 2.16e0  &  35.75   & $c/20$  & 2.23e0 &  5.61   & 2.07e0   & 149   \cr
    \hline
    $1$&1e-3  &$c/(5M)$& 5.78e-4&  1019.15   &$c/(30n)$& 5.90e-4&  5604.80& 5.68e-4  & 99   \cr
       &1e-2  &$c/(5M)$& 1.14e-2 &246.40    &$c/(30n)$& 1.15e-2& 1356.87  & 1.12e-2  &24     \cr
       &5e-2  &$c/(5M)$& 6.47e-2 & 112.20   &$c/(30n)$& 6.48e-2&  613.41 & 6.19e-2  & 11    \cr
    \hline
    $2$&1e-3  &$c/(10M)$& 7.57e-5& 474.10   &$c/(50n)$& 1.32e-4&  2441.85& 6.82e-5  & 23    \cr
       &1e-2  &$c/(10M)$& 1.80e-3 &229.90    &$c/(50n)$& 1.92e-3&1047.03   & 1.47e-3  &10     \cr
       &5e-2  &$c/(10M)$& 2.32e-2& 156.75   &$c/(50n)$& 2.35e-2&  708.72 & 1.61e-2  & 6    \cr
      \hline
    $4$&1e-3  &$c/(10M)$& 2.51e-5&  250.80  &$c/(60n)$& 1.03e-4&  2212.26& 1.30e-5  & 10    \cr
       &1e-2  &$c/(10M)$& 1.14e-3 &170.50    &$c/(60n)$& 1.29e-3&941.19   & 6.42e-4  &6      \cr
       &5e-2  &$c/(10M)$& 2.23e-2& 138.05   &$c/(60n)$& 2.25e-2&  746.67 & 8.58e-3  & 3     \cr
      \bottomrule
    \end{tabular}
    \end{threeparttable}
\end{table}

Now we examine more closely the convergence behaviour of the SVRG iterates, and compare
it with that of SGD and LM. For all these three examples and all $\nu$ values, both SVRG
and SGD can achieve an accuracy comparable with that by LM, {thereby achieving
the order optimality of these methods,} when the step size $c_0$
for SVRG and SGD is taken to be of order $\mathcal{O}(M^{-1})$ and $\mathcal{O}(n^{-1})$, respectively. This
observation agrees well with the analysis in Theorem \ref{thm:main}. Generally, the larger the index $\nu$ is,
the smaller the value $c_0$ should be taken in order to achieve the optimal rate. This can also
be seen partly from the constant $2^{2\nu}c_\nu$ in the error bound in Theorem
\ref{thm:main}. Next we discuss the computational complexity. For all three examples,
SVRG takes fewer epochs to reach the optimal error than SGD for a large index $\nu$, and
LM requires fewest iterations among the three methods. For
small $\nu$, SVRG stops earlier than LM, and can be faster than SGD for suitably chosen $c_0$
(see, e.g., the case $\nu=0$ in Table \ref{tab:phil}). {These empirical observations agree with the
fact that SVRG hybridizes SGD and LM.} Since in practice the index $\nu$ is rarely known,
SVRG is an excellent choice, due to its low sensitivity with respect to $\nu$.

\begin{table}[htp!]
  \centering\small
  \begin{threeparttable}
  \caption{Comparison between SVRG (with $M=100$), SGD and LM for \texttt{s-shaw}.\label{tab:shaw}}
    \begin{tabular}{ccccccccccccccc}
    \toprule
    \multicolumn{2}{c}{Method}&
    \multicolumn{3}{c}{SVRG}&\multicolumn{3}{c}{SGD}&\multicolumn{2}{c}{LM}\cr
    \cmidrule(lr){3-5} \cmidrule(lr){6-8}
    \cmidrule(lr){9-10}
    $\nu$& $\epsilon$ &$c_0$&$e_{\rm svrg}$&$k_{\rm svrg}$&$c_0$&$e_{\rm sgd}$&$k_{\rm sgd}$&$e_{\rm lm}$&$k_{\rm lm}$\\
    \midrule
    $0$&1e-3 &$  c  $& 2.81e-1  & 30246.15   &$ c   $& 2.81e-1  & 2704.92   & 2.81e-1  & 760983 \cr
       &1e-2 & $c$ & 6.92e-1 &503.25  &$ c   $& 7.08e-1   & 42.42  &6.67e-1    & 12385   \cr
       &5e-2 & $c$ & 3.01e0 & 139.15 &$ c   $& 3.91e0   &10.59   &2.91e0    & 3392   \cr
    \hline
    $1$&1e-3 &$ c/M $& 6.80e-5  & 579.15     &$c/(2n) $& 7.05e-5 & 1047.60   & 5.95e-5  & 144    \cr
       &1e-2 & $c/M$ & 5.35e-3  &222.75      &$c/(2n) $& 5.42e-3 & 394.00   & 5.21e-3  &54
       \cr
       &5e-2 & $c/M$ & 1.50e-1  &  148.50    &$c/(2n)$& 1.50e-1 & 271.00      & 1.47e-1  & 36     \cr
    \hline
    $2$&1e-3 & $c/(2M)$ & 6.94e-5  & 434.50  &$c/(20n)$ & 7.08e-5 & 4147.00   & 6.36e-5  & 50    \cr
       &1e-2 & $c/(2M)$ & 5.80e-3  &246.95    &$c/(20n) $& 5.80e-3 & 2242.50     & 5.71e-3  &30       \cr
        &5e-2 & $c/(2M)$ & 7.84e-2  & 52.80   &$c/(20n) $& 7.79e-2  & 480.80     & 7.08e-2  & 5      \cr
    \hline
    $4$&1e-3 & $c/(4M)$ & 3.83e-5  & 184.25  &$c/(30n)$&5.79e-5 & 1966.38   & 3.13e-5  & 9     \cr
       &1e-2 & $c/(4M)$ &1.96e-3  &121.55     &$c/(30n)$&1.99e-3 &828.45     & 1.01e-3  &4       \cr
       &5e-2 & $c/(4M)$ &3.61e-2  & 95.15    &$c/(30n)$&3.61e-2 & 645.75    & 6.45e-3  & 1      \cr
    \bottomrule
    \end{tabular}
    \end{threeparttable}
\end{table}

\begin{figure}[hbt!]
\centering
\setlength{\tabcolsep}{4pt}
\begin{tabular}{ccc}
\includegraphics[width=0.31\textwidth,trim={1.5cm 0 0.5cm 0.5cm}]{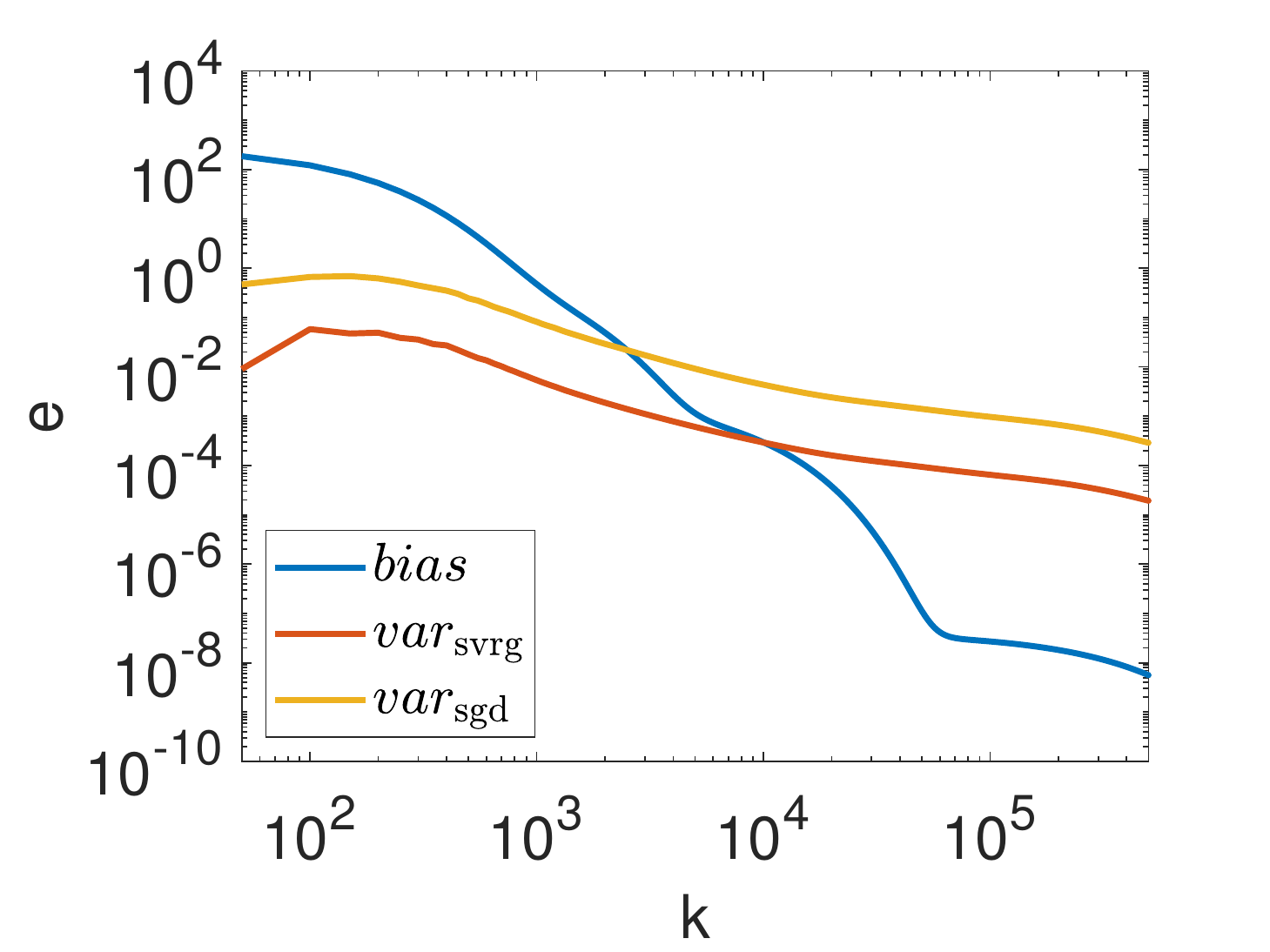}&
\includegraphics[width=0.31\textwidth,trim={1.5cm 0 0.5cm 0.5cm}]{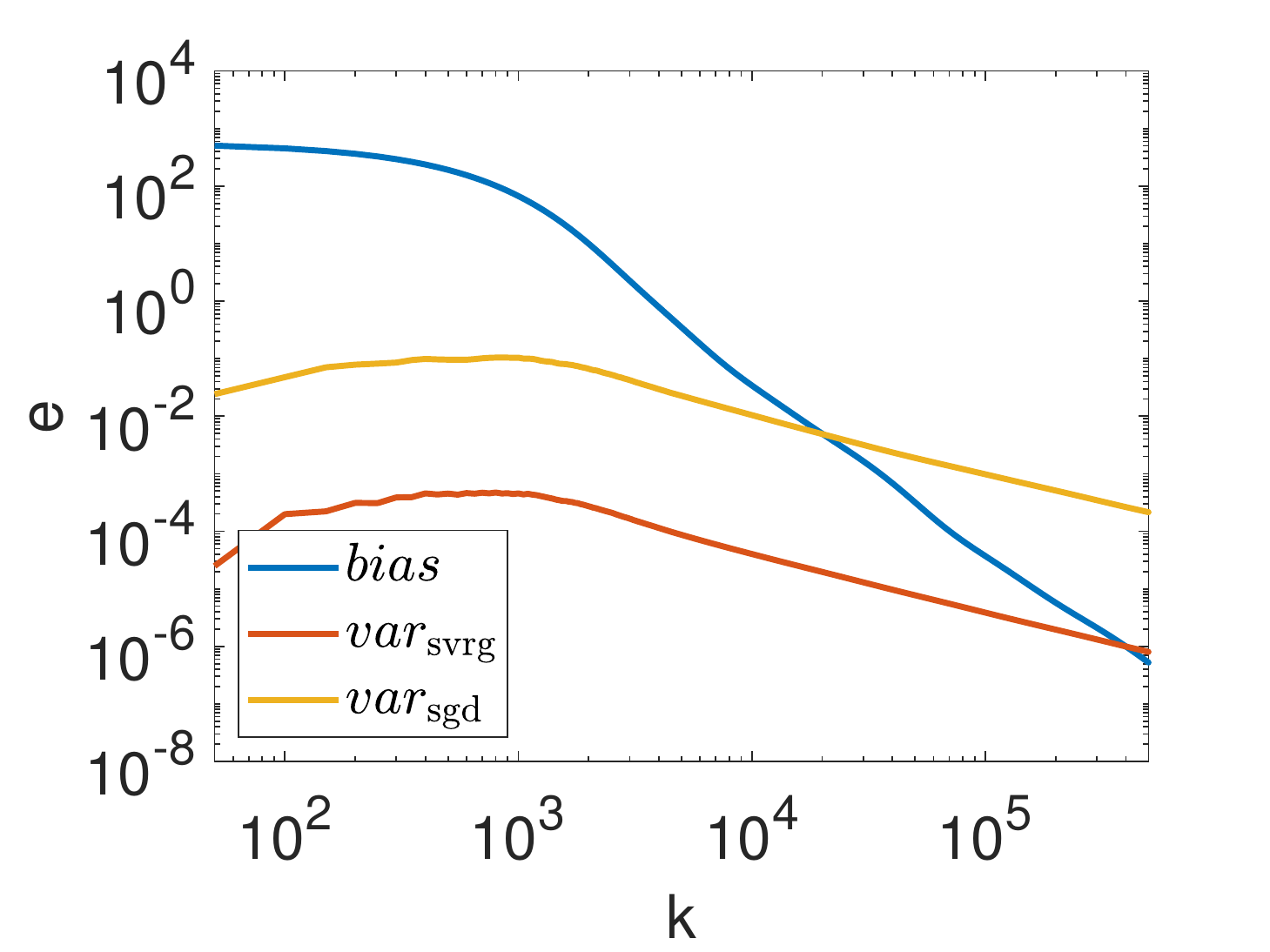}&
\includegraphics[width=0.31\textwidth,trim={1.5cm 0 0.5cm 0.5cm}]{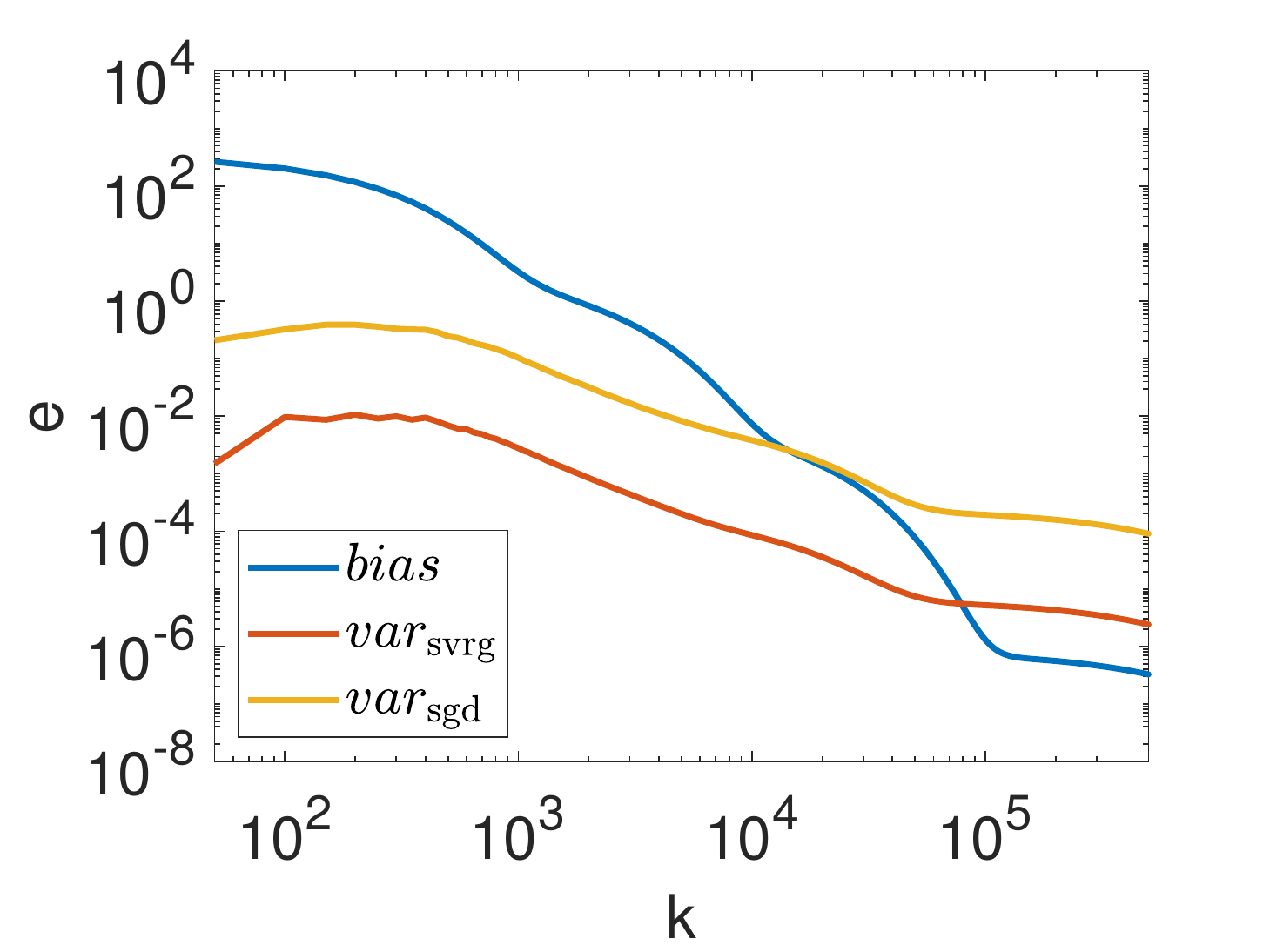}\\
\includegraphics[width=0.31\textwidth,trim={1.5cm 0 0.5cm 0.5cm}]{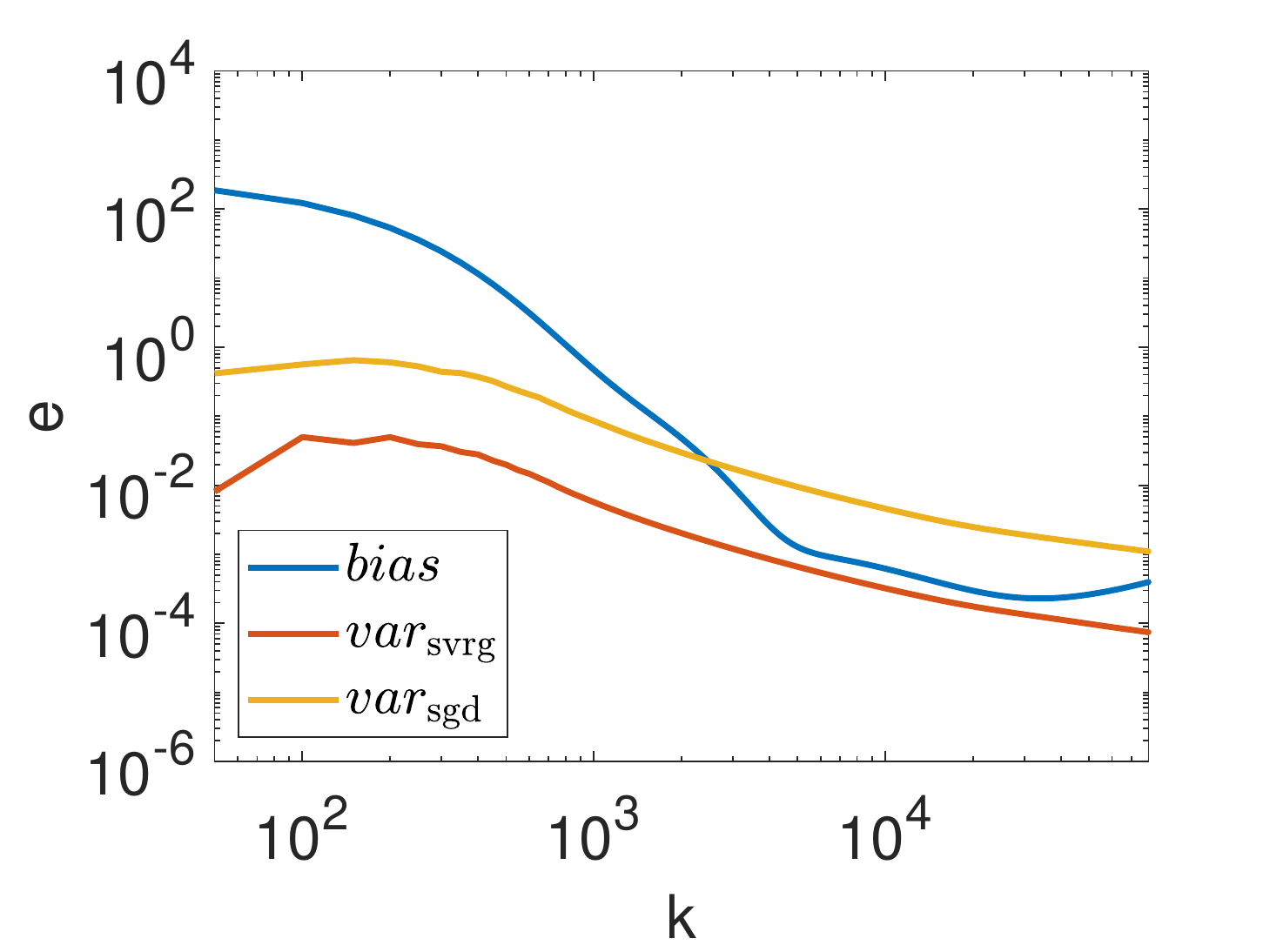}&
\includegraphics[width=0.31\textwidth,trim={1.5cm 0 0.5cm 0.5cm}]{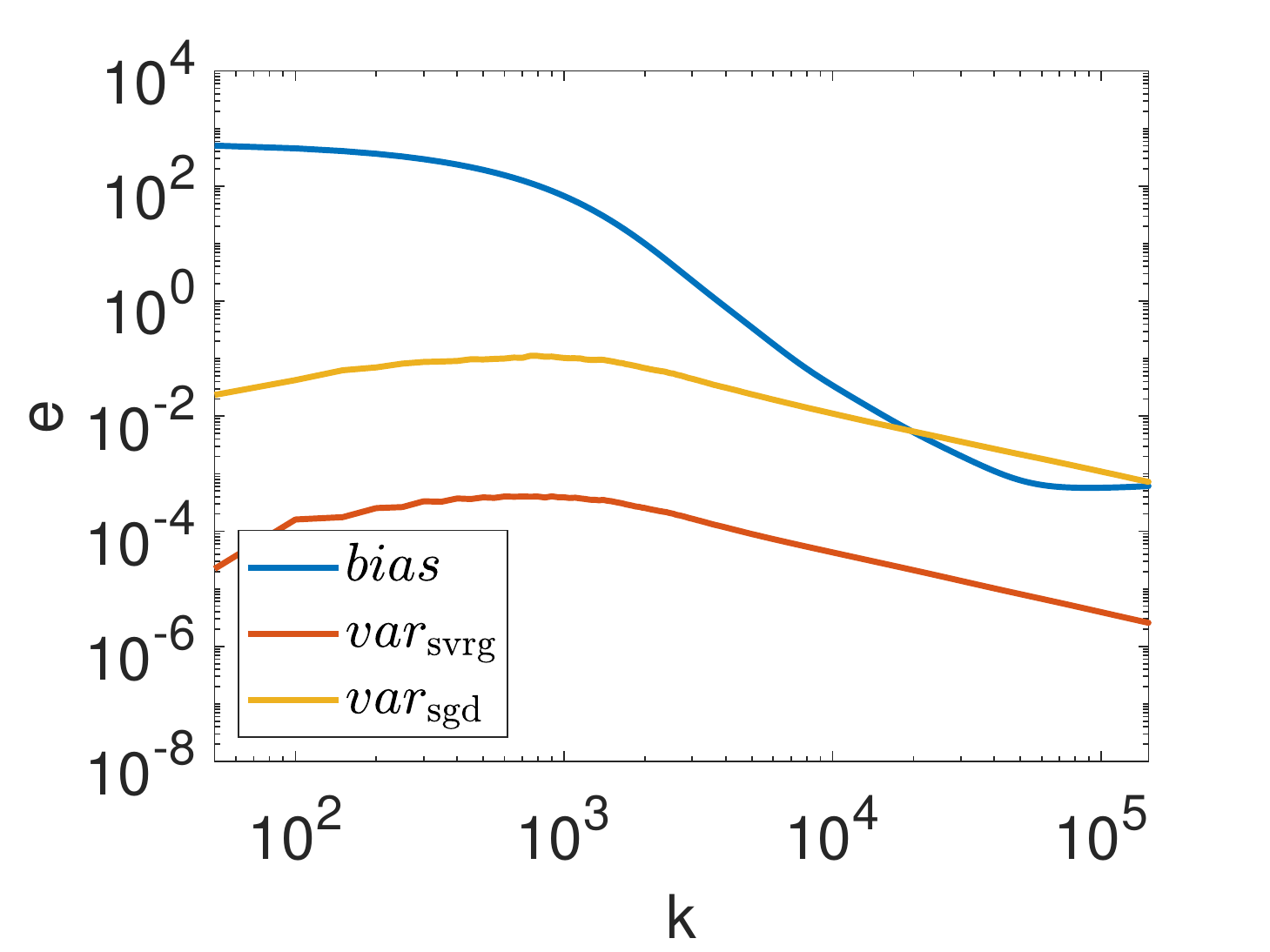}&
\includegraphics[width=0.31\textwidth,trim={1.5cm 0 0.5cm 0.5cm}]{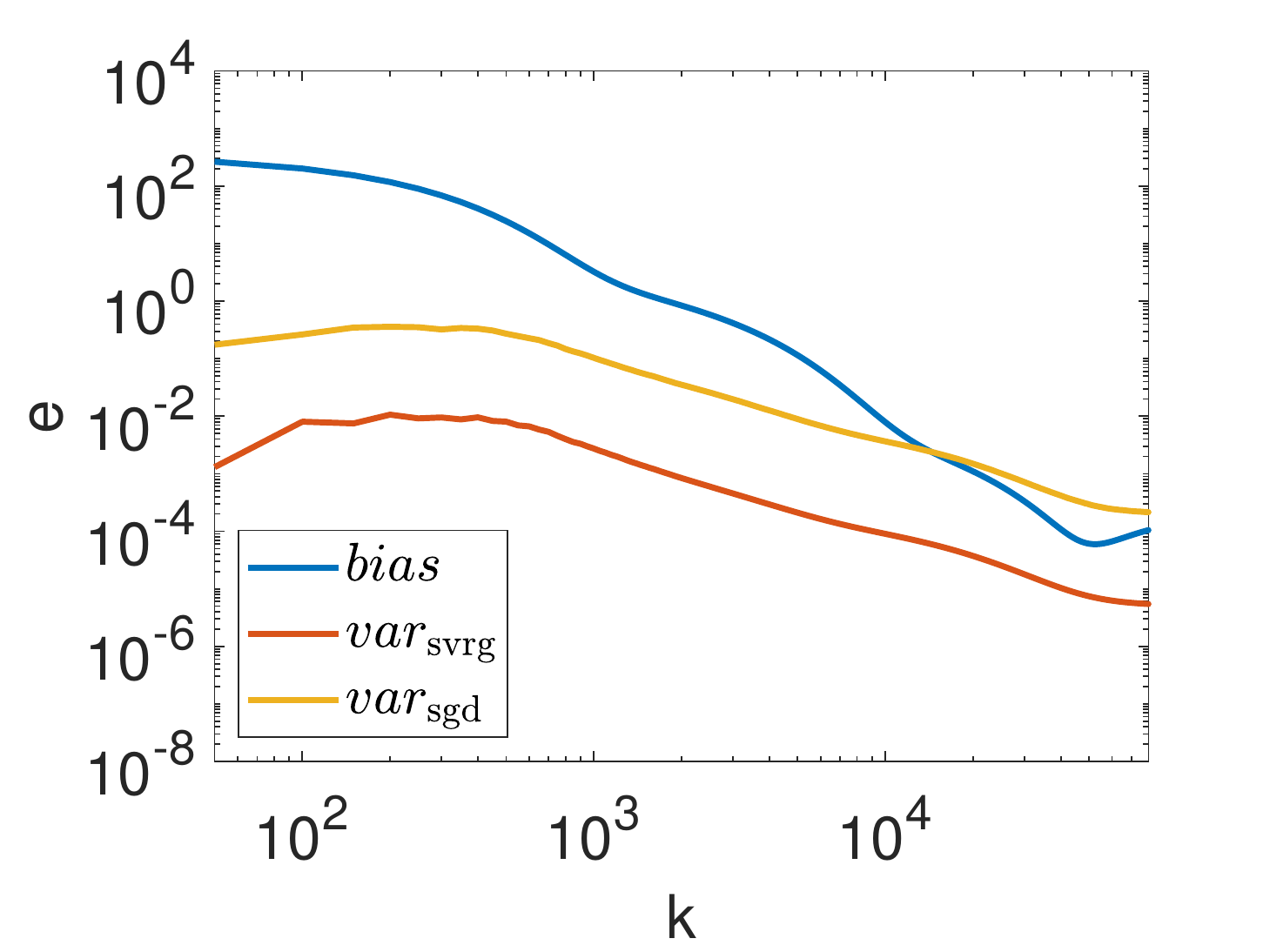}\\
\includegraphics[width=0.31\textwidth,trim={1.5cm 0 0.5cm 0.5cm}]{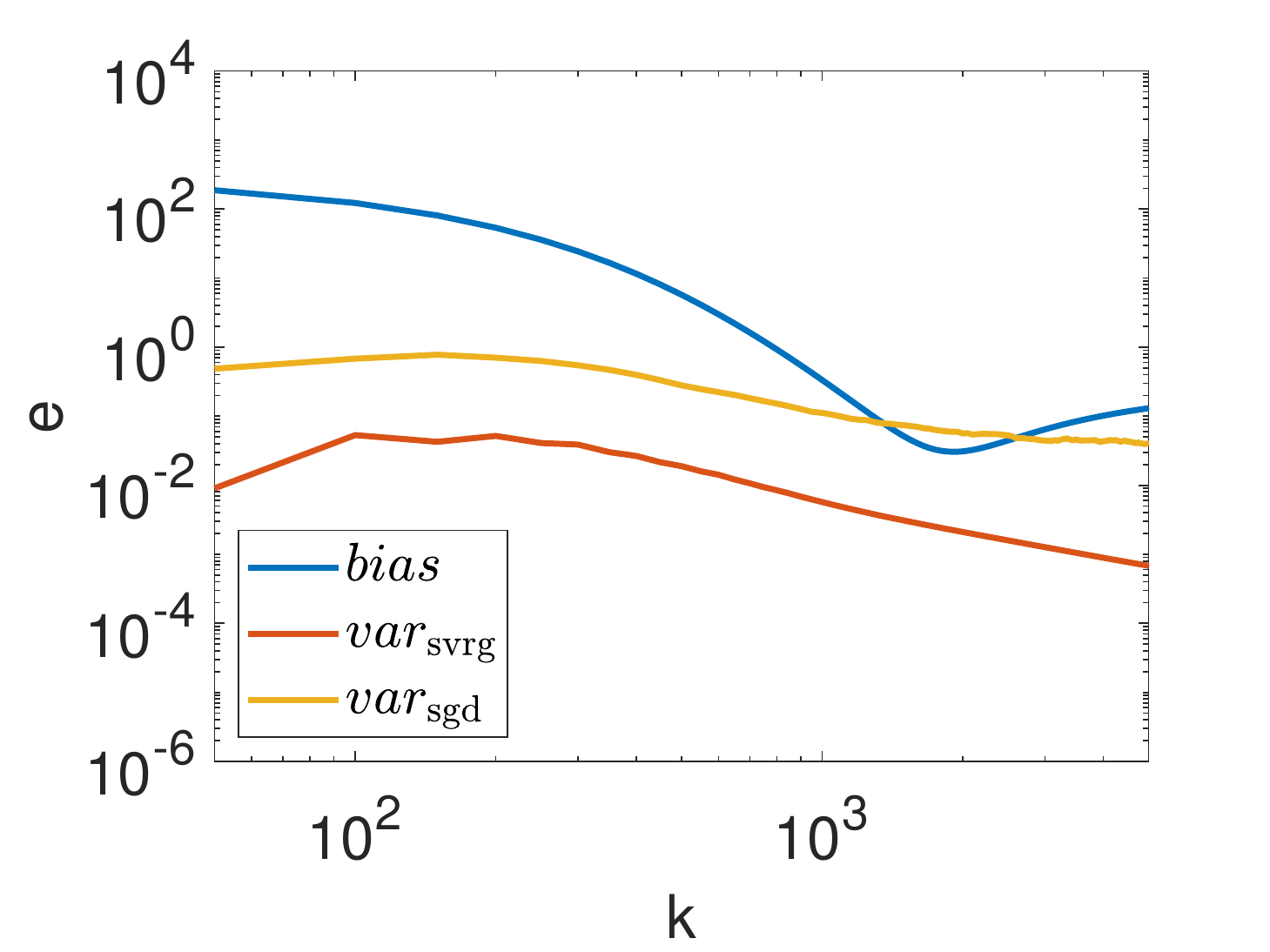}&
\includegraphics[width=0.31\textwidth,trim={1.5cm 0 0.5cm 0.5cm}]{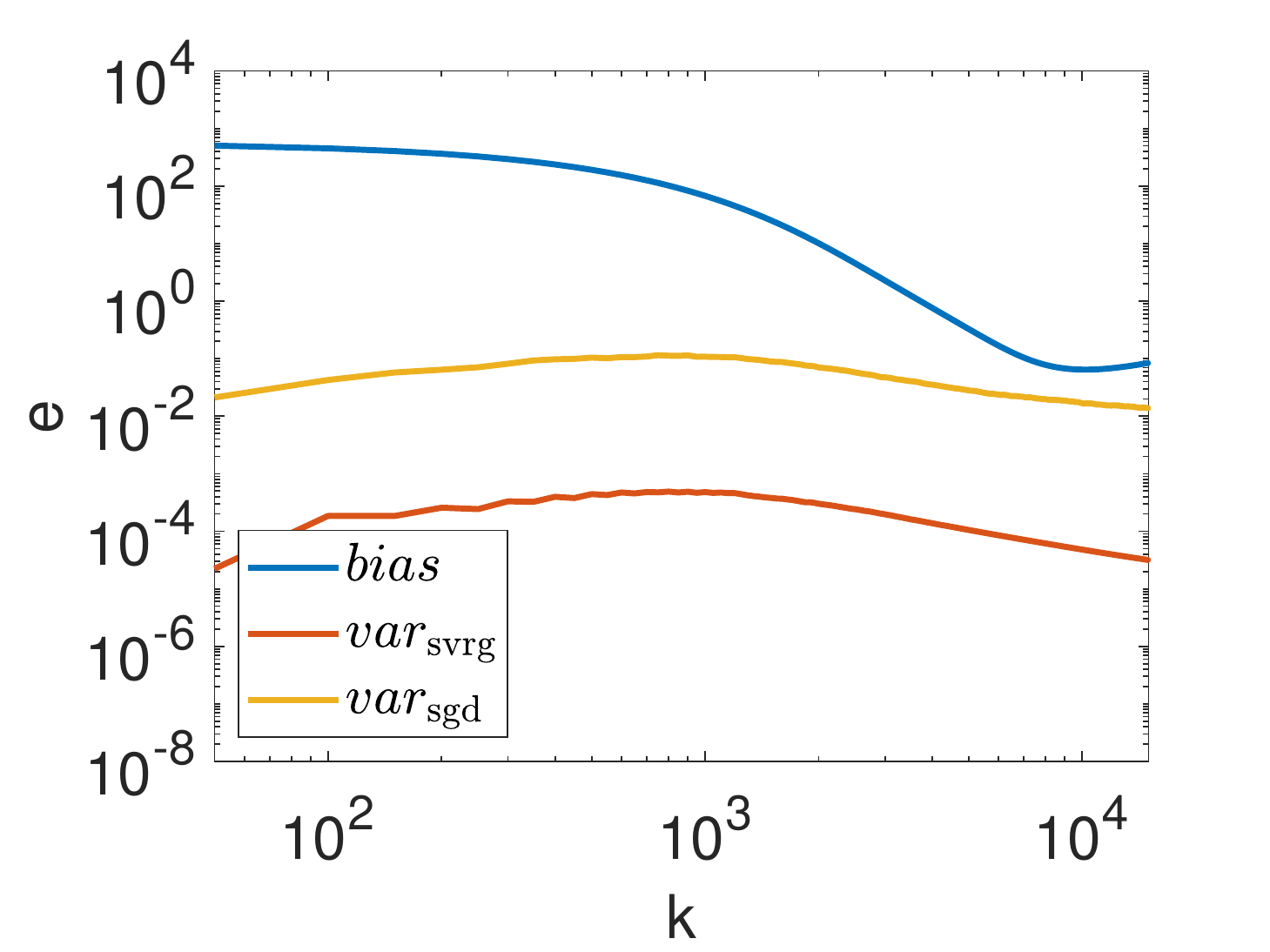}&
\includegraphics[width=0.31\textwidth,trim={1.5cm 0 0.5cm 0.5cm}]{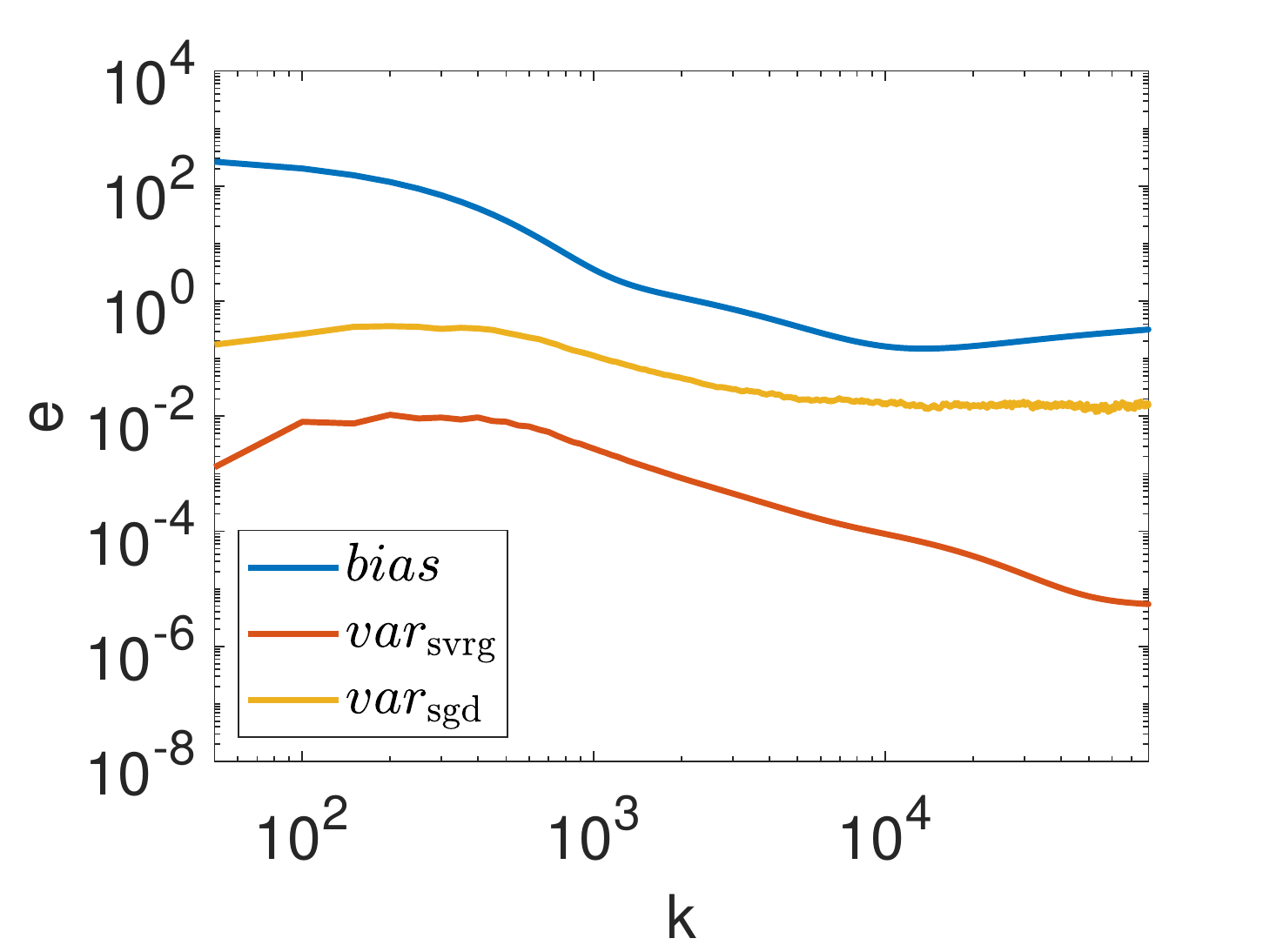}\\
\texttt{s-phillips}& \texttt{s-gravity} & \texttt{s-shaw}
\end{tabular}
\caption{The convergence of the bias or variance with generic term $e$ versus iteration number for the examples with
$\nu=1$. The rows from top to bottom rows are for $\epsilon=0$,
$\epsilon=$1e-3 and $\epsilon$=5e-2, respectively.\label{fig:err}}
\end{figure}

To verify the analysis in Section \ref{sec:comp}, we examine the bias $bias=\|\E[x_k^\delta] - x^\dag\|^2=\|\E[\hat{x}_k^\delta] - x^\dag\|^2$,
and the variances $var_{\rm svrg}=\E[\|x_k^\delta-\E[x_k^\delta]\|^2]$ and $var_{\rm sgd}=\E[\|\hat{x}_k^\delta-\E[\hat{x}_k^\delta]\|^2]$.
The numerical results are shown in Fig. \ref{fig:err}, for the examples with $\nu=1$, with the
step size $c_0$ for SVRG used for both methods. Although not presented, we note that any other suitable $c_0$ under condition
\eqref{cond:compare} leads to nearly identical observations. Note that the iteration index $k$ in the figures refers to the exact number
of iterations (not counted in epoch), to facilitate the comparison of the convergence behaviour. For both exact and noisy data, when
the iteration number $k$ is fixed, the
SVRG variance $var_{\rm svrg}$ is always orders of magnitude smaller than the SGD variance $var_{\rm sgd}$, which is fully in line with Theorem \ref{thm:svrg-sgd}. This shows clearly the role of the variance reduction effect, which in particular allows using larger step size.
Note that the frequency $M=100$ is selected by the condition \eqref{cond:sat_svrg} for optimal accuracy, but actually does not satisfy
condition \eqref{cond:compare}. Nonetheless, we still observe the assertion in Theorem \ref{thm:svrg-sgd}.

Further, in the experiments, $bias$ (which is equal to the error $e_{\rm lm}$ of Landweber method) is
always much larger than the SVRG variance $var_{\rm svrg}$ (of similar magnitude during {a few iterations before stopping}),
and thus the variance has little influence on the optimal accuracy, especially for noisy data. In contrast, the SGD variance
$var_{\rm sgd}$ dominates the error sometimes and causes the undesirable saturation phenomenon. These observations also
agree with Theorem \ref{thm:main}, which states that the saturation of SVRG does not exist by choosing suitable
frequency $M$ and initial step size $c_0$. They also confirm the theoretical prediction in Remark \ref{rem:sat_0}, i.e., the
condition for the optimality of SVRG is weaker than that of SGD, partly concurring with Theorem \ref{thm:svrg-sgd}.
{These empirical observations show clearly the beneficial effect of incorporating variance reduction
into stochastic iterative methods from the perspective of regularization theory.}

\subsection{Influence of $M$}\label{sec:M}
{SVRG involves one free parameter, the frequency $M$ of evaluating the full gradient.
Clearly, the parameter $M$ will influence the overall computational efficiency of SVRG: ideally
one would like to make it as large as possible, but a too large $M$ would bring
too little variance reduction into SGD iteration. The theoretical analysis in this work
indicates that SVRG can achieve optimal convergence rates when $M\geq \mathcal{O}(n^\frac12)$
(cf. Remark \ref{rem:con_optimal}), and that $M\leq\mathcal{O}(n^{\frac12})$ is
sufficient for ensuring the SVRG variance smaller than SGD variance (cf. Remark \ref{cond:M_comp}). Nonetheless,
a complete theoretical analysis of the influence of the frequency $M$ on the
performance of SVRG is still unknown. To gain insight, we present the numerical results
for \texttt{s-phillips} with noisy data by SVRG with different $M$ ranging from
$0.1n$ to $5n$ in Table \ref{tab:phil_M}. Note that the choices $2n$ and $5n$ were
recommended for convex and nonconvex optimization problems, respectively \cite{JohnsonZhang:2013}.
The numerical results indicate that SVRG with all these frequencies can actually
achieve an accuracy comparable with that by the Landweber
method when the constant step size is chosen suitably. In general, a larger $M$ requires
smaller step sizes in order to maintain the optimal convergence rate, agreeing well with
the theoretical analysis in Section \ref{sec:conv}. Interestingly, the overall computational
complexity for these different $M$ does not vary too much. Thus, the choice of $M$
within a certain range actually has little impact on the performance of SVRG. Although not
presented, the same observations can be drawn from the numerical results for the examples \texttt{s-shaw}
and \texttt{s-gravity}.}

\begin{table}[htp!]
  \centering\small
  \begin{threeparttable}
  \caption{SVRG with different $M$ for \texttt{s-phillips}.\label{tab:phil_M}}
    \begin{tabular}{cccccccccccc}
    \toprule
    \multicolumn{2}{c}{}&
    \multicolumn{3}{c}{$\nu=0$}& \multicolumn{3}{c}{$\nu=2$}&\\
    \cmidrule(lr){3-5} \cmidrule(lr){6-8}
    $M$& $\epsilon$ &$c_0$&$e$&$k$&$c_0$&$e$&$k$&\\
    \midrule
    $0.1n$&1e-3 & $5c/M$  &1.67e-2& 4134.35   & $c/(2M)$ & 7.16e-5 & 155.10  \cr
    &1e-2 & $5c/M$  & 1.31e-1 &180.95  & $c/(2M)$  & 1.07e-3 &68.75  \cr
       &5e-2 & $5c/M$  & 5.42e-1 &96.80 & $c/(2M)$   & 2.90e-2 & 46.75 \cr
    \hline
    $0.5n$&1e-3 & $5c/M$  &1.66e-2& 5650.35   & $c/(2M)$ & 4.18e-5 & 204.30  \cr
    &1e-2 & $5c/M$  & 1.31e-1 &125.70  & $c/(2M)$  & 9.90e-4 &93.30  \cr
       &5e-2 & $5c/M$  & 5.40e-1 &66.15 & $c/(2M)$   & 2.90e-2 & 63.75 \cr
    \hline
    $n$&1e-3 & $10c/M$  &1.67e-2& 3757.40   & $c/M$ & 5.83e-5 & 139.50  \cr
    &1e-2 & $10c/M$  & 1.29e-1 &163.80
  & $c/M$  & 1.04e-3 &62.20  \cr
       &5e-2 & $10c/M$  & 5.38e-1 &87.40 & $c/M$   & 2.92e-2 & 42.50 \cr
    \hline
    $2n$&1e-3 & $15c/M$  &1.67e-2& 3781.35   & $1.5c/M$ & 7.63e-5 & 144.38  \cr
    &1e-2 & $15c/M$  & 1.30e-1 &164.70  & $1.5c/M$  & 1.08e-3 &62.25  \cr
       &5e-2 & $15c/M$  & 5.39e-1 &87.08 & $1.5c/M$   & 2.93e-2 & 42.53 \cr
    \hline
    $5n$&1e-3 & $25c/M$  &1.66e-2& 4519.86   & $2c/M$ & 7.33e-5 & 214.32  \cr
    &1e-2 & $25c/M$  & 1.29e-1 &197.28  & $2c/M$  & 1.05e-3 &93.60  \cr
       &5e-2 & $25c/M$  & 5.40e-1 &104.64 & $2c/M$   & 2.90e-2 & 63.84 \cr
    \bottomrule
    \end{tabular}
    \end{threeparttable}
\end{table}

\subsection{On Assumption \ref{ass:stepsize}(iii)}\label{sec:condition}

\begin{figure}[hbt!]
\centering
  \setlength{\tabcolsep}{4pt}
\begin{tabular}{ccc}
\includegraphics[width=0.31\textwidth,trim={1.5cm 0 0.5cm 0.5cm}]{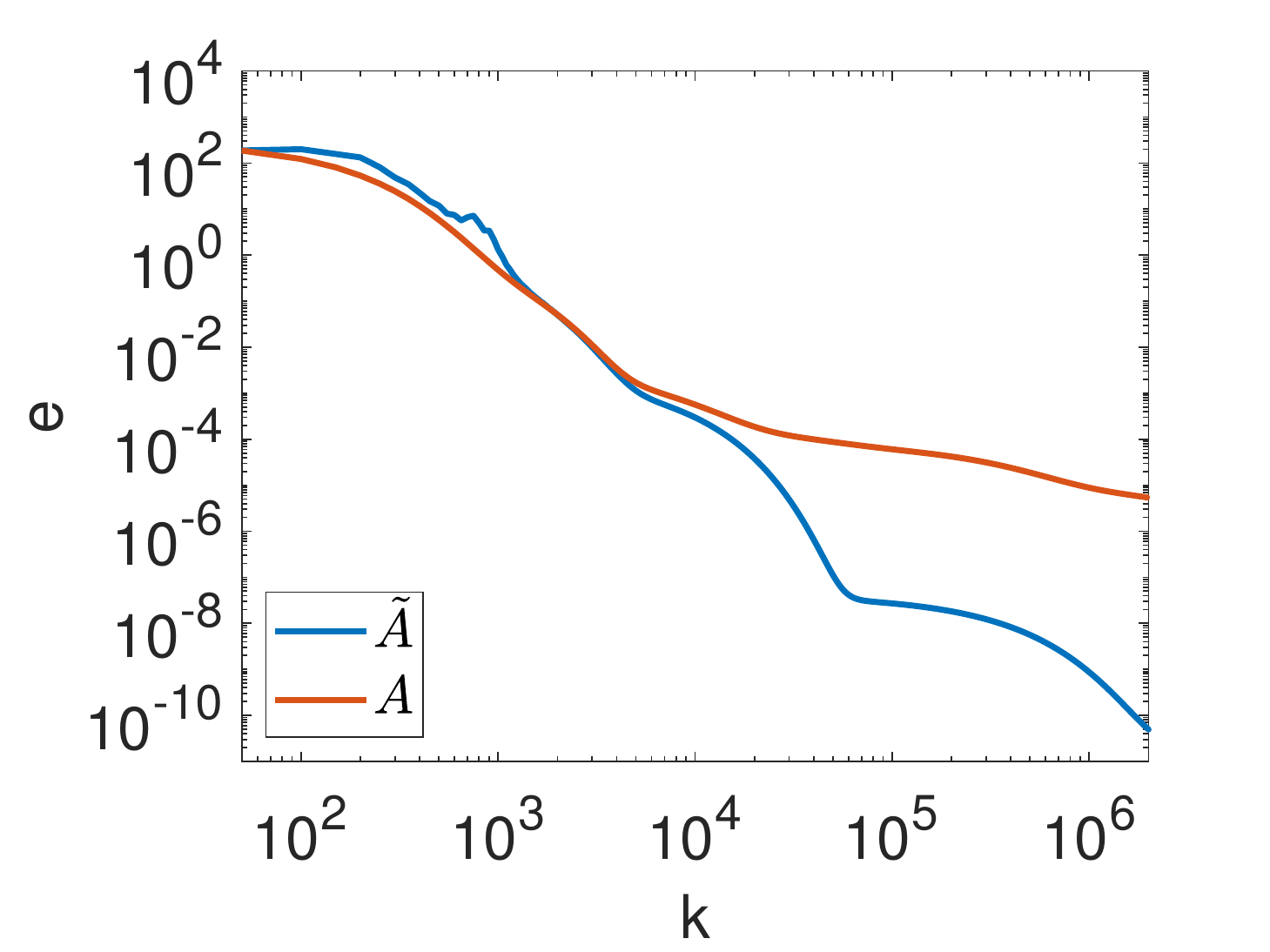}&
\includegraphics[width=0.31\textwidth,trim={1.5cm 0 0.5cm 0.5cm}]{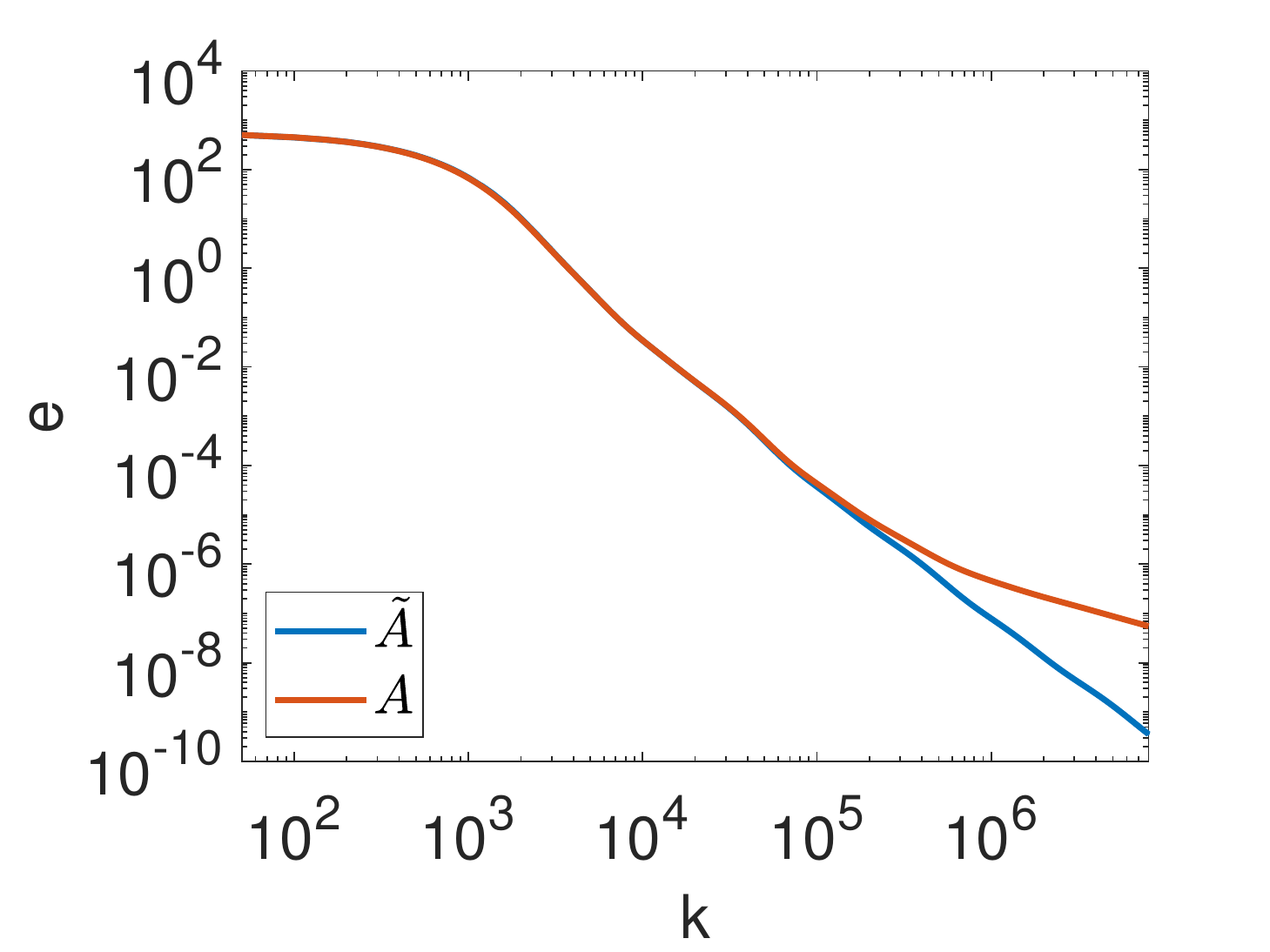}&
\includegraphics[width=0.31\textwidth,trim={1.5cm 0 0.5cm 0.5cm}]{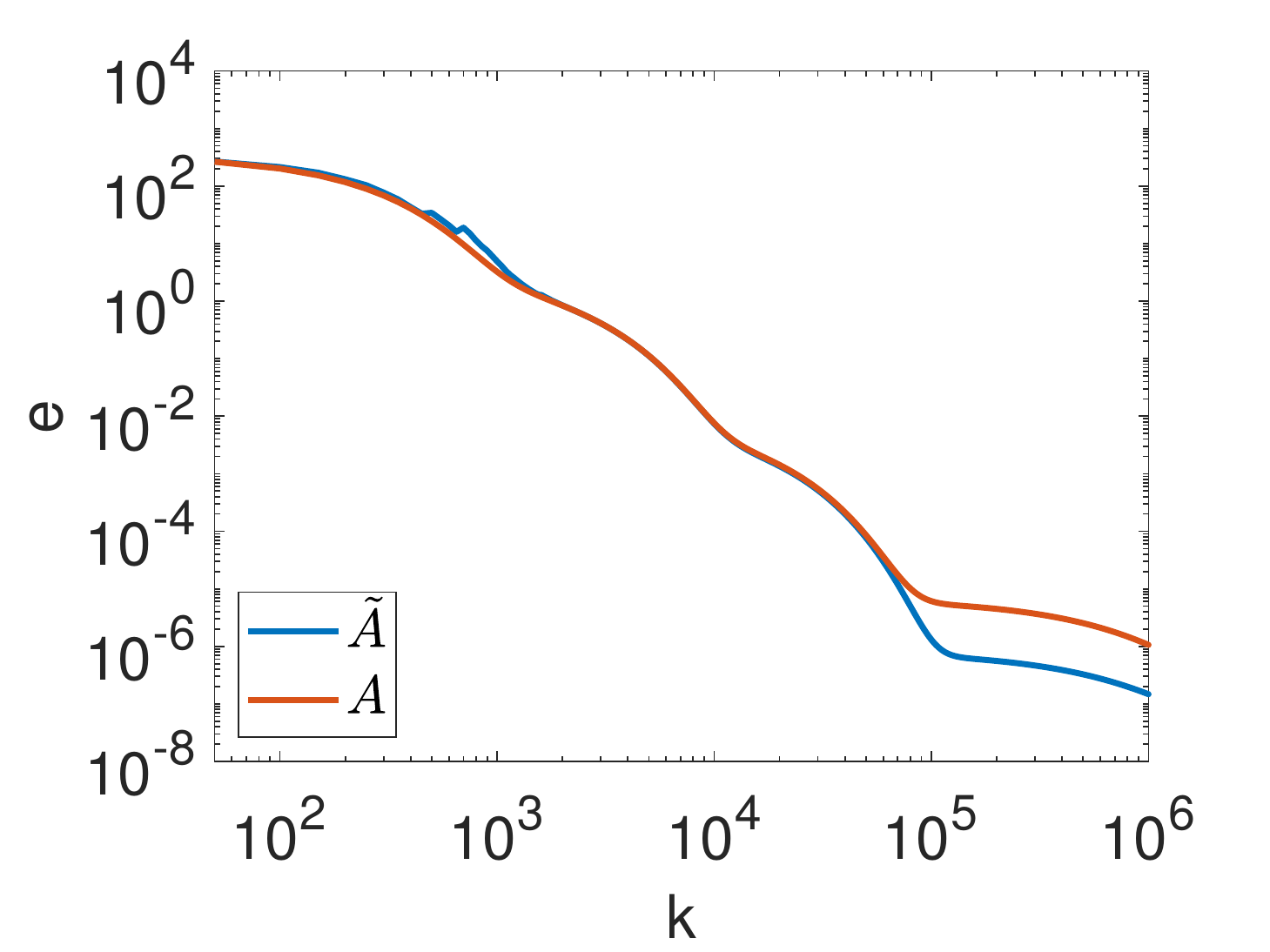}\\
\includegraphics[width=0.31\textwidth,trim={1.5cm 0 0.5cm 0.5cm}]{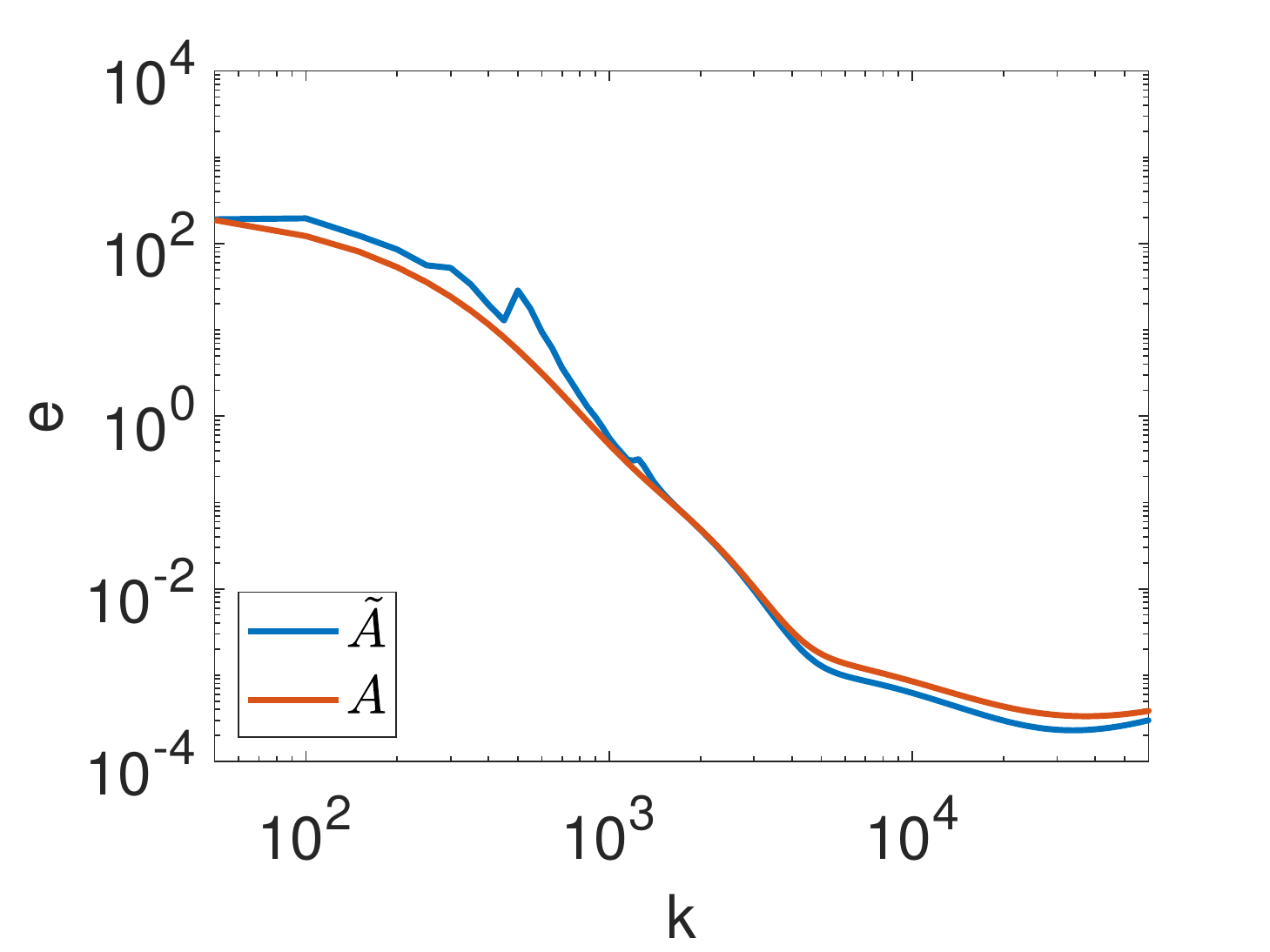}&
\includegraphics[width=0.31\textwidth,trim={1.5cm 0 0.5cm 0.5cm}]{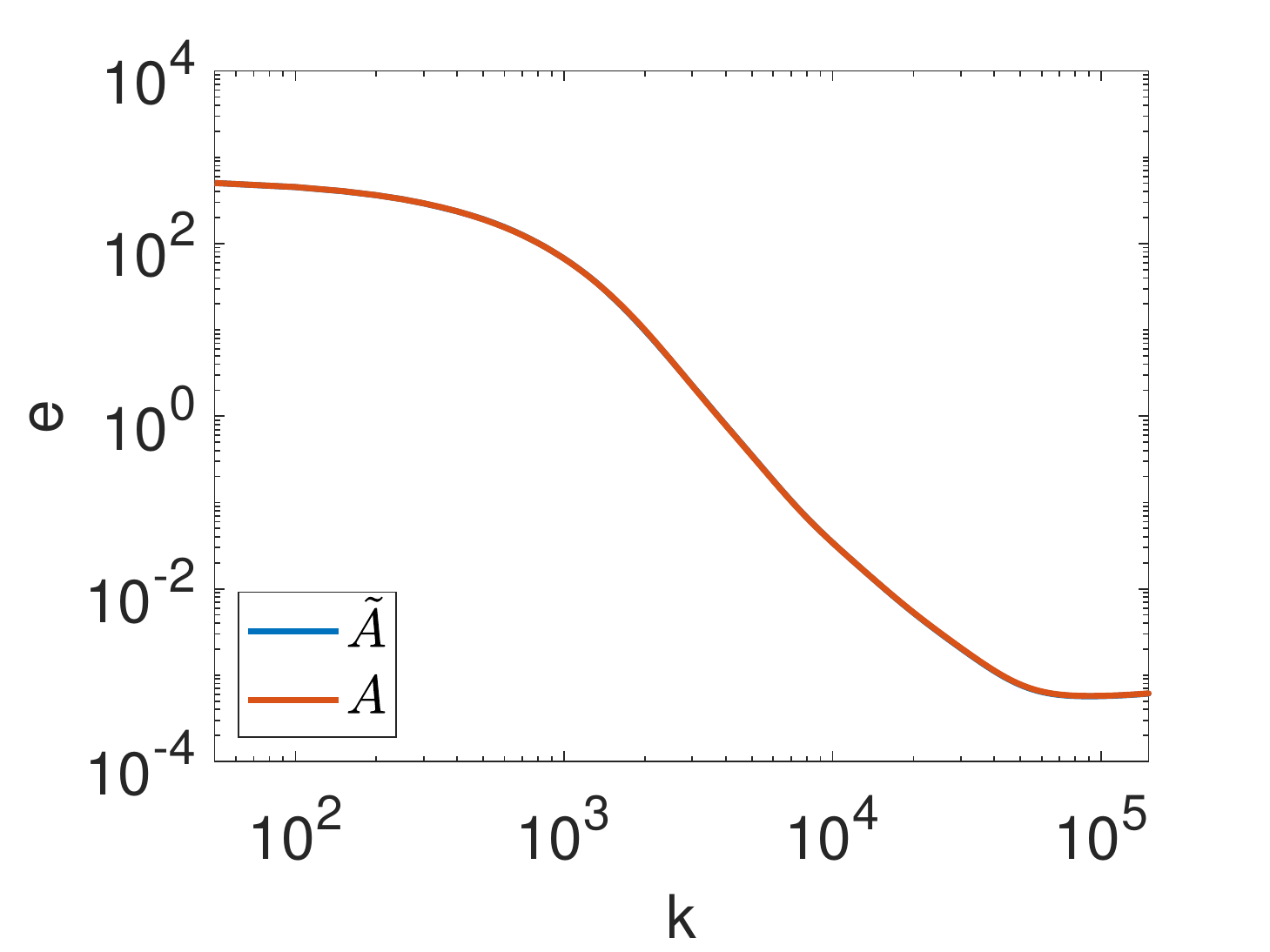}&
\includegraphics[width=0.31\textwidth,trim={1.5cm 0 0.5cm 0.5cm}]{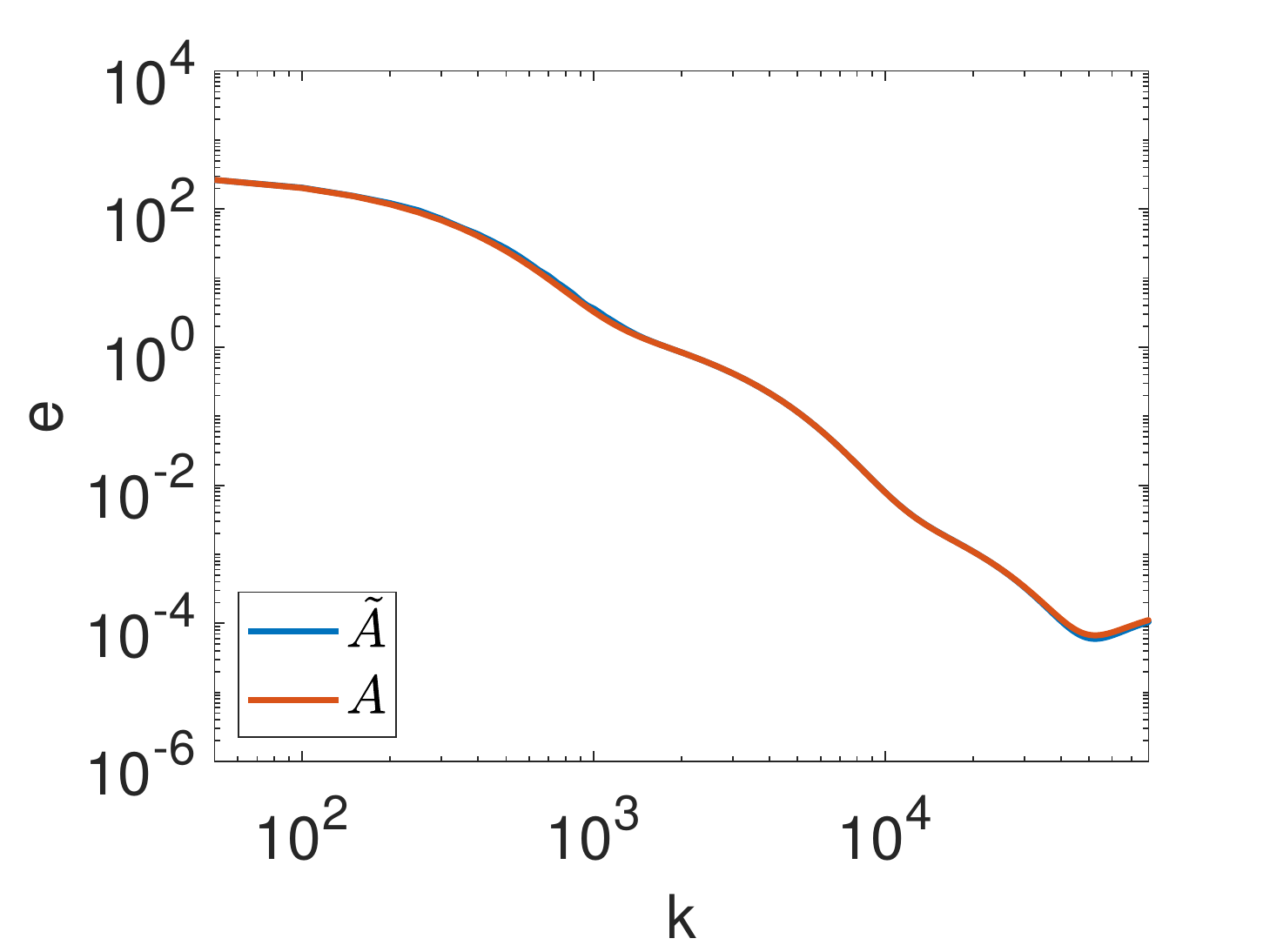}\\
\includegraphics[width=0.31\textwidth,trim={1.5cm 0 0.5cm 0.5cm}]{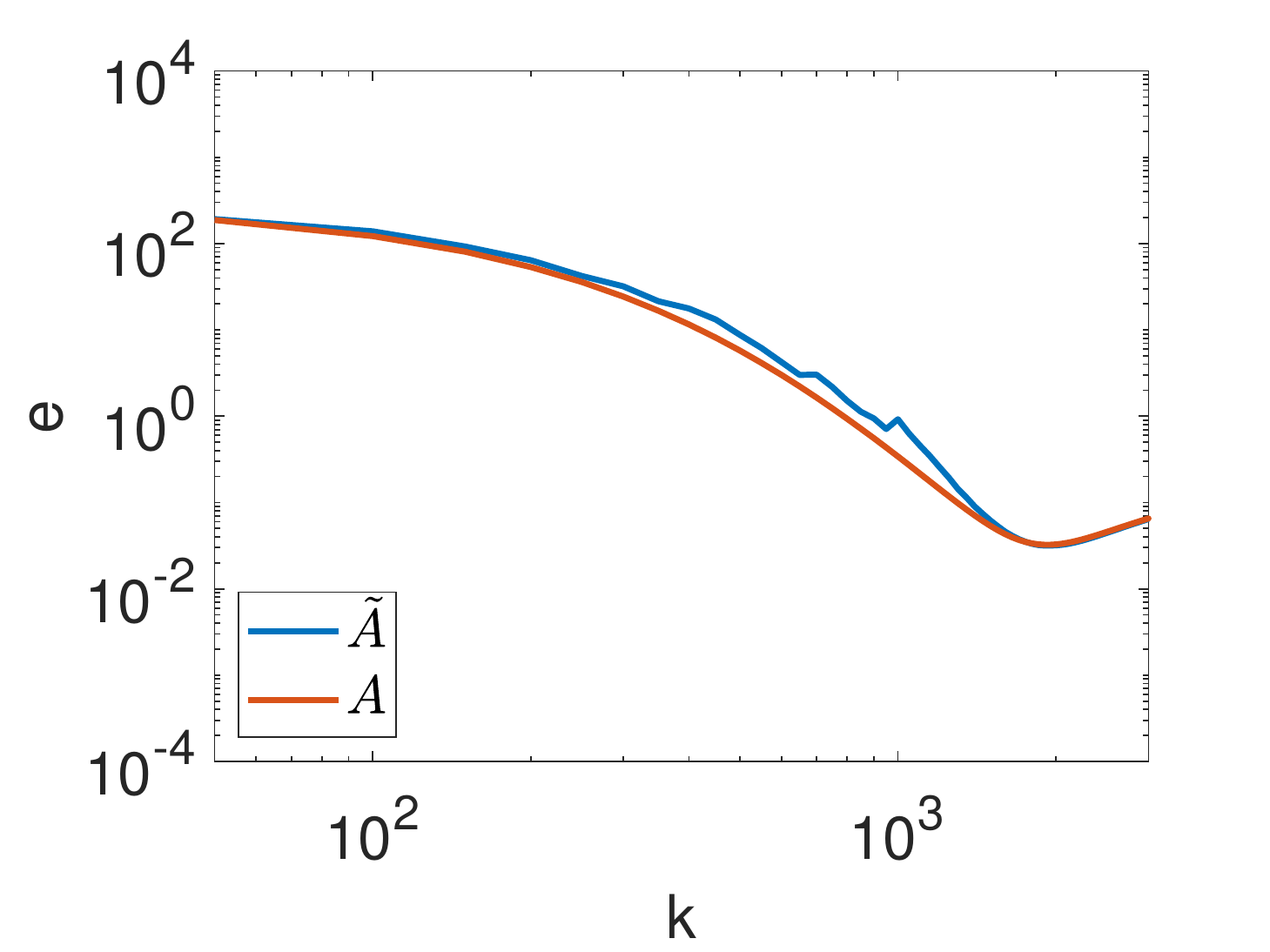}&
\includegraphics[width=0.31\textwidth,trim={1.5cm 0 0.5cm 0.5cm}]{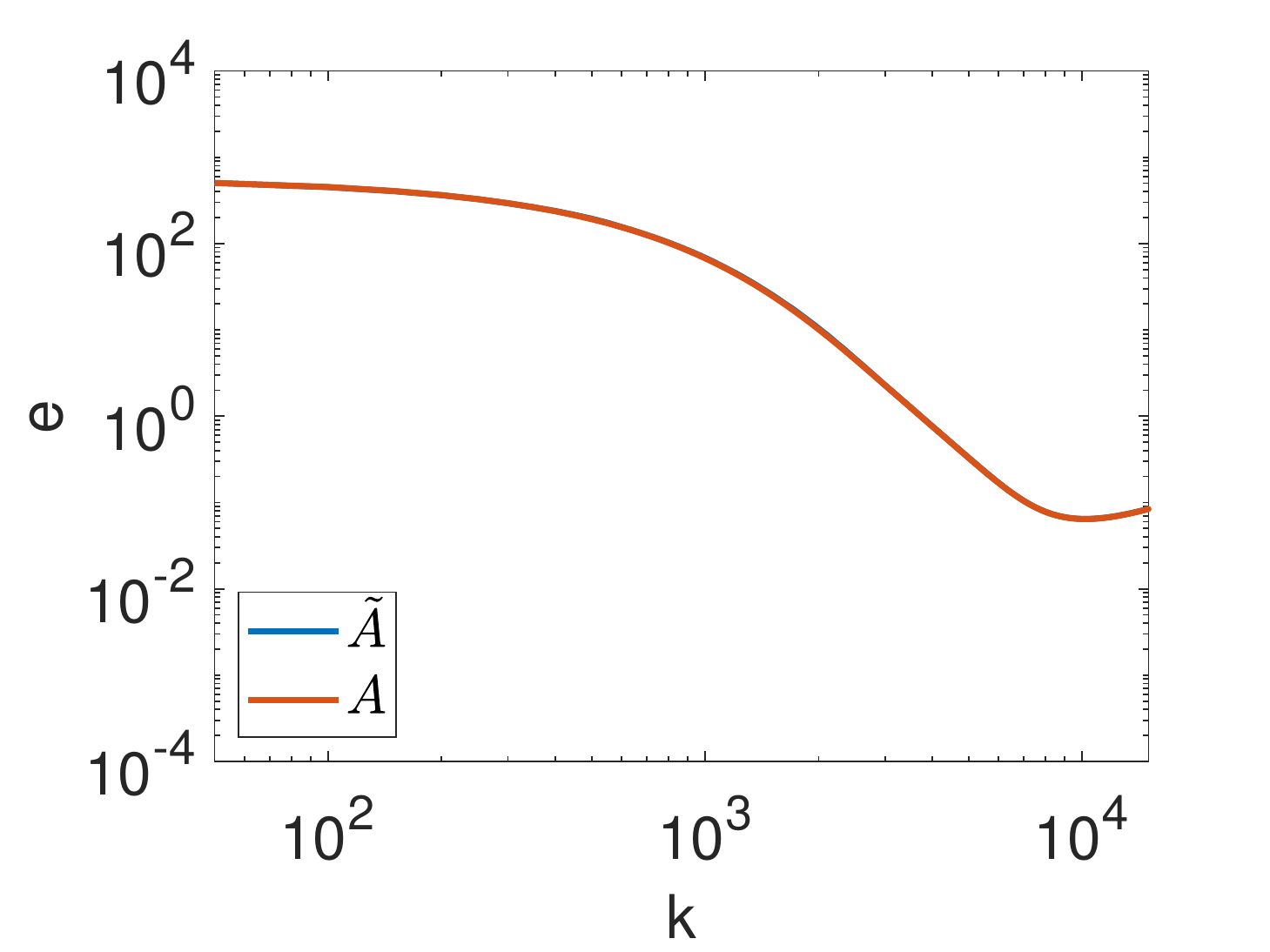}&
\includegraphics[width=0.31\textwidth,trim={1.5cm 0 0.5cm 0.5cm}]{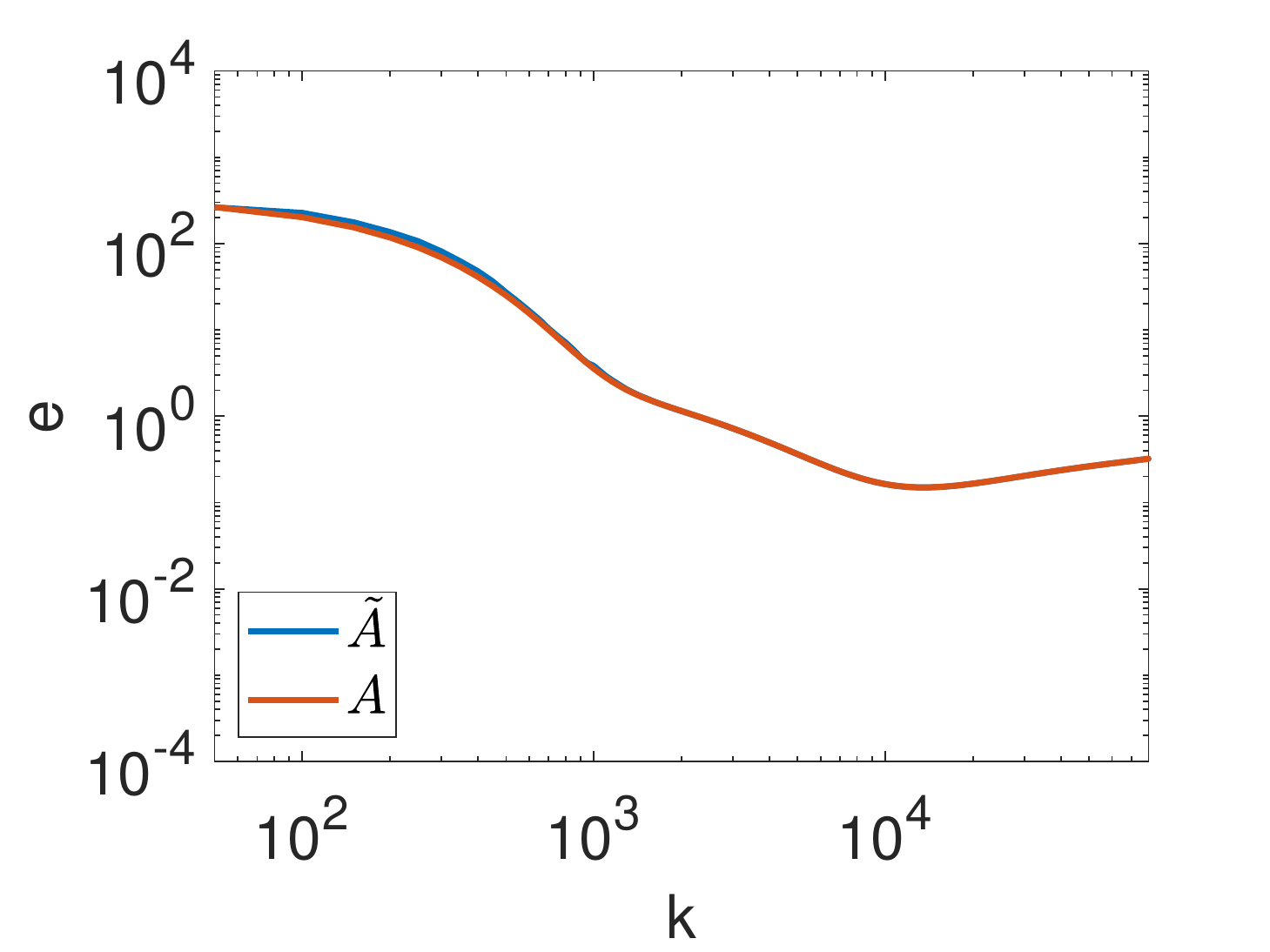}\\
\texttt{s-phillips}& \texttt{s-gravity} & \texttt{s-shaw}
\end{tabular}
\caption{The convergence of the error $e$ versus iteration number for the examples with
$\nu=1$, computed using $A$ and $\tilde A$. The rows from top to bottom rows are for $\epsilon=0$,
$\epsilon=$1e-3 and $\epsilon$=5e-2, respectively.\label{fig:err_UA}}
\end{figure}

Assumption \ref{ass:stepsize}(iii) is crucial to the analysis in Sections \ref{sec:conv} and \ref{sec:comp}.
{It is natural to ask whether the assumption is actually necessary. We examine the issue numerically
as follows.} Let $A=U\Sigma V^t$ be the SVD  of $A$, and $\tilde A$ by $\tilde A = U^tA$, and then replace $A$ in
\eqref{eqn:lininv} by $\tilde A$ and $y^\delta$ by $\tilde y^\delta = U^ty^\delta$. Then preconditioned
system $\tilde A x = \tilde y^\delta$ satisfies Assumption \ref{ass:stepsize}(iii).
The numerical results for \texttt{s-phillips} are shown in Table \ref{tab:phil_UA}, and
the trajectories of $e_k^\delta$ for the examples with $\nu=1$ in Fig. \ref{fig:err_UA}.
It is observed that for noisy data, the SVRG results for $A$ and $\tilde A$ are nearly
identical with each other in terms of the accuracy, stopping index, and convergence
trajectory. For exact data (cf. the top row of Fig. \ref{fig:err_UA}), the trajectories
overlap up to a certain point {around 1e-3 for \texttt{s-phillips} and 1e-5 for
\texttt{s-gravity} and \texttt{s-shaw},} which can be further decreased
by choosing smaller $c_0$. These observations {resemble closely the empirical observations for SGD, see, especially Fig. 4.3 of \cite{JinZhouZou:2021}.}
Thus, Assumption \ref{ass:stepsize}(iii) is probably due to a limitation of the proof technique,
and there might be  alternative proof strategies that circumvent the restriction.

\begin{table}[htp!]
  \centering\small
  \begin{threeparttable}
  \caption{Comparison between SVRG (with $M=100$) for \texttt{s-phillips} with $A$ and $\tilde A$.\label{tab:phil_UA}}
    \begin{tabular}{cccccccccccc}
    \toprule
    \multicolumn{3}{c}{Method}&
    \multicolumn{2}{c}{ SVRG with $A$}& \multicolumn{2}{c}{ SVRG with $\tilde A$}&\\
    \cmidrule(lr){4-5} \cmidrule(lr){6-7}
    $\nu$& $\epsilon$ &$c_0$&$e$&$k$&$e$&$k$&\\
    \midrule
    $0$&1e-3 & $5c/M$  &1.67e-2& 4134.35   & 1.65e-2 & 4129.40  \cr
    &1e-2 & $5c/M$  & 1.31e-1 &180.95    & 1.28e-1 &176.55  \cr
       &5e-2 & $5c/M$  & 5.42e-1 &96.80    & 5.36e-1 & 96.25 \cr
    \hline
    $1$&1e-3 & $1.5c/M$   & 3.31e-4 &  430.65  & 2.29e-4 & 372.35\cr
       &1e-2 & $1.5c/M$   &  5.96e-3 &41.25 & 5.32e-3 &40.70  \cr
       &5e-2 & $1.5c/M$   & 3.22e-2 & 21.45  & 3.17e-2 & 20.90 \cr
    \hline
    $2$&1e-3 &$c/(2M)$& 7.16e-5 &  155.10 & 3.49e-5&  148.50\cr
       &1e-2 &$c/(2M)$& 1.07e-3 &68.75   & 9.77e-4 &68.75  \cr
       &5e-2 &$c/(2M)$& 2.90e-2 &  46.75  & 2.89e-2& 46.75 \cr
    \hline
    $4$&1e-3 &$c/(5M)$& 3.05e-5 &  202.95 & 2.46e-5& 201.30 \cr
       &1e-2 &$c/(5M)$& 2.41e-3 &142.45    &2.41e-3 &142.45  \cr
       &5e-2 &$c/(5M)$& 5.20e-2 &  110.00  & 5.21e-2& 110.00 \cr
    \bottomrule
    \end{tabular}
    \end{threeparttable}
\end{table}

\appendix

\section{Technical proofs}

In this appendix, we collect the proofs of several technical estimates.

\subsection{Proof of Lemma \ref{lem:kernel}}

The proof relies on spectral decomposition. Let $\mathrm{Sp}(B)$ be the spectrum of $B$.
Then by direct computation, we have
\begin{align*}
c_0^s\|B^s M_0^{KM}\|=&c_0^s\sup_{\lambda\in{\rm Sp}(B)}|\lambda^s (1-c_0\lambda)^{KM}|\leq \sup_{a\in[0,1]}a^s(1-a)^{KM}.
\end{align*}
Let $g(a)=a^s(1-a)^{KM}$. Then $g'(a)=\big(s(1-a)-KMa\big)a^{s-1}(1-a)^{KM-1}$, so that
$g(a)$ achieves its maximum over the interval $[0,1]$ at $a^*=s(s+KM)^{-1}$. Consequently,
\begin{align*}
c_0^s\|B^s M_0^{KM}\|\leq &g(a^*)=(\tfrac{KM}{s+KM})^{s+KM}s^s (KM)^{-s}\leq s^s M^{-s} K^{-s}.
\end{align*}
This shows the second estimate. Similarly,
\begin{align*}
c_0^{-t}\|B^{-t} (I-M_0^{KM})\|=&\sup_{\lambda\in{\rm Sp}(B)}|(c_0\lambda)^{-t} (1-(1-c_0\lambda)^{KM})|
\leq\sup_{a\in[0,1]}a^{-t} (1-(1-a)^{KM}).
\end{align*}
Note that for any $a\in[0,1]$, there holds {$1-(1-a)^{KM}\leq1$, and}
$\min_{t\in[0,1]}(aKM)^t=\min(aKM,1)$,
since $(aKM)^t$ is monotone with respect to $t$. Let $h(a):=aKM-(1-(1-a)^{KM})$ which is increasing
over $[0,1]$, that implies $h(a)\geq h(0)=0$. Thus
$${1-(1-a)^{KM}\leq \min(aKM,1)\leq (aKM)^t}.$$
This shows the first estimate and completes the proof of the lemma.

\subsection{Proof of Proposition \ref{prop:bias-var}}
To prove Proposition \ref{prop:bias-var}, we
first give a representation of the (epochwise) SVRG iterate $x_{KM}^\delta$.
\begin{lemma}\label{lem:recursion}
The following recursion holds for any $K\geq 0$,
\begin{align}
e_{(K+1)M}^\delta =
(M_0^M-L_KB) e_{KM}^\delta+\Big(c_0\sum_{i=0}^{M-1}M_0^i+L_K\Big)\zeta, \label{eqn:SVRG_K}
\end{align}
where the random matrix $L_K$ is given by
\begin{equation}\label{eqn:Lj}
  L_K=c_0\sum_{i=1}^{M-1}H_{KM+i}(I-M_0^i)B^{-1}.
\end{equation}
\end{lemma}
\begin{proof}
Note that the SVRG iterate $x_{k+1}^\delta$, $k=0,1,\ldots$, can be rewritten as
\begin{align*}
  x_{k+1}^\delta &= x_k^\delta - c_0 \big((a_{i_k} ,e_k^\delta-e_{k_M}^\delta)a_{i_k}+Be_{k_M}^\delta-\zeta\big)\nonumber\\
                 &= x_k^\delta - c_0 a_{i_k}a_{i_k}^t (e_k^\delta-e_{k_M}^\delta) - c_0(Be_{k_M}^\delta-\zeta).
\end{align*}
Using the definitions of $P_k$ and $N_k$, the error $e_{k}^\delta\equiv x_{k}^\delta -x^\dag$ of the SVRG iterate $x_k^\delta$ satisfies
\begin{align}
  e_{k+1}^\delta &= (I-c_0a_{i_k}a_{i_k}^t)e_k^\delta + c_0(a_{i_k}a_{i_k}^t-B)e_{k_M}^\delta + c_0\zeta
                 = P_ke_k^\delta - c_0N_ke_{k_M}^\delta  + c_0\zeta.\label{eqn:svrg-it-0}
\end{align}
For any $K\geq 0$, it follows from \eqref{eqn:svrg-it-0} and direct computation that
\begin{equation}
e_{KM+1}^\delta = P_{KM}e_{KM}^\delta-c_0N_{KM}e_{KM}^\delta+c_0\zeta=M_0e_{KM}^\delta+c_0\zeta.\label{eqn:svrg-it-KM}
\end{equation}
Meanwhile, setting $k=(K+1)M-1$ in the recursion \eqref{eqn:svrg-it-0}, then repeatedly applying the
recursion \eqref{eqn:svrg-it-0} and using the definitions of the matrices $G_k$ and $H_k$ lead to
\begin{align*}
  e_{(K+1)M}^\delta =& P_{(K+1)M-1}e_{(K+1)M-1}^\delta - c_0N_{(K+1)M-1} e_{KM}^\delta  + c_0\zeta\nonumber\\
                    =& G_{(K+1)M-2} e_{(K+1)M-2}^\delta
                      - c_0(P_{(K+1)M-1}N_{(K+1)M-2} \\
                     & + N_{(K+1)M-1})e_{KM}^\delta  + c_0(P_{(K+1)M-1}+I)\zeta \nonumber\\
                    = & ...
                    = G_{KM+1}e_{KM+1}^\delta-c_0\sum_{i=1}^{M-1}H_{KM+i}e_{KM}^\delta+c_0\sum_{i=2}^{M}G_{KM+i}\zeta. 
\end{align*}
This identity and \eqref{eqn:svrg-it-KM} imply that for any $K\geq 0$,
\begin{align*}
e_{(K+1)M}^\delta =&\Big(G_{KM+1}M_0-c_0\sum_{i=1}^{M-1}H_{KM+i}\Big)e_{KM}^\delta+\Big(c_0\sum_{i=1}^{M}G_{KM+i}\Big)\zeta. 
\end{align*}
Next we simplify the two terms in the {brackets} using the identity
\eqref{eqn:G}. It follows directly from \eqref{eqn:G} that
\begin{align*}
  G_{KM+1}M_0-c_0\sum_{i=1}^{M-1}H_{KM+i} = M_0^M-c_0\sum_{i=1}^{M-1}H_{KM+i}(I-M_0^i).
\end{align*}
Similarly, by the identity \eqref{eqn:G}, we deduce
\begin{align*}
c_0\sum_{i=1}^{M}G_{KM+i}
=&c_0I+c_0\sum_{i=1}^{M-1}\Big(M_0^{M-i}+c_0\sum_{j=0}^{M-i-1}H_{KM+i+j}M_0^j\Big)\\
=&c_0\sum_{i=1}^{M}M_0^{M-i}+c_0^2\sum_{i=1}^{M-1}\sum_{j=0}^{M-i-1}H_{KM+i+j}M_0^j\\
=&c_0\sum_{i=0}^{M-1}M_0^{i}+c_0^2\sum_{i=1}^{M-1}H_{KM+i}\Big(\sum_{j=0}^{i-1}M_0^j\Big)\\
=&c_0\sum_{i=0}^{M-1}M_0^{i}+c_0\sum_{i=1}^{M-1}H_{KM+i}(I-M_0^i)B^{-1},
\end{align*}
where the last line follows from the identity
\begin{equation}\label{eqn:sum-M}
c_0\sum_{i=0}^{j-1}M_0^i=(I-M_0^{j})B^{-1}
\end{equation}
Combining the preceding identities completes the proof of the lemma.
\end{proof}

Now we can give the proof of Proposition \ref{prop:bias-var}.
\begin{proof}
By the definitions of the matrices $N_i$ and $G_{i+1}$, they are independent. Thus, there hold
$$\E[H_i]=\E[G_{i+1}]\E[N_i]=0\quad \mbox{and}\quad \E[L_j]=0.$$
Then by Lemma \ref{lem:recursion}, we have
\begin{align*}
\E[e_{(K+1)M}^\delta]=M_0^M \E[e_{KM}^\delta]+c_0\sum_{i=0}^{M-1}M_0^i \zeta.
\end{align*}
Repeatedly applying this identity gives
\begin{align*}
 &\E[e_{(K+1)M}^\delta]=M_0^M \Big(M_0^M \E[e_{(K-1)M}^\delta]+c_0\sum_{i=0}^{M-1}M_0^i \zeta\Big)+c_0\sum_{i=0}^{M-1}M_0^i \zeta\\
=&M_0^{2M} \E[e_{(K-1)M}^\delta]+c_0\sum_{i=0}^{2M-1}M_0^i \zeta
=\cdots=M_0^{(K+1)M} e_0^\delta+c_0\sum_{i=0}^{(K+1)M-1}M_0^i \zeta.
\end{align*}
This and the identity \eqref{eqn:sum-M} show the expression for $\E[e_{KM}^\delta]$.
Let $z_{K}:=e_{KM}^\delta-\E[e_{KM}^\delta].$
Then for any $K\geq 0$, it follows from Lemma \ref{lem:recursion} that
\begin{align*}
z_{K+1}=&M_0^M z_K+ R_K, \quad\mbox{with }R_K:=L_K(\zeta-B e_{KM}^\delta),
\end{align*}
and $z_0=0.$ Repeatedly applying the recursion directly gives
\begin{align*}
z_{K+1}=M_0^{(K+1)M} z_{0}+\sum_{j=0}^{K}M_0^{jM}  R_{K-j}=\sum_{j=0}^K M_0^{(K-j)M} R_{j}.
\end{align*}
This completes the proof of the proposition.
\end{proof}

\subsection{Proof of Proposition \ref{prop:bias-var_sgd}}
The following recursion is direct from the definition of SGD iteration in \eqref{eqn:SGD}
\begin{equation*}
  \hat e_{k+1}^\delta = (I-c_0a_{i_k}a_{i_k}^t)\hat e_k^\delta + c_0\xi_{i_k}a_{i_k} = P_k \hat e_k^\delta + c_0\zeta_k.
\end{equation*}
Repeatedly applying the recursion and using the identity \eqref{eqn:G} (and its proof) yield that for any $K\geq 0$,
\begin{align*}
\hat{e}_{(K+1)M}^\delta =&G_{KM+1}P_{KM}\hat{e}_{KM}^\delta+c_0\sum_{i=0}^{M-1}G_{KM+i+1}\zeta_{KM+i}\\
=&\Big(M_0^M+c_0\sum_{i=0}^{M-1}H_{KM+i}M_0^i\Big)\hat{e}_{KM}^\delta+c_0\zeta_{(K+1)M-1}\\
&+c_0\sum_{i=1}^{M-1}\Big(M_0^{M-i}+c_0\sum_{t=0}^{M-i-1}H_{KM+i+t}M_0^t\Big)\zeta_{KM+i-1}.
\end{align*}
Since $\E[H_{KM+i}]=0$, for $i=0,\ldots,M-1$,
and $H_{KM+i+t}$, $t\geq0$, and $\zeta_{KM+i-1}$ are independent, by the identity \eqref{eqn:sum-M},
\begin{align*}
\E[\hat{e}_{(K+1)M}^\delta]=&M_0^M\E[\hat{e}_{KM}^\delta]+c_0\sum_{i=0}^{M-1}M_0^{i}\zeta
=M_0^{(K+1)M}\hat{e}_{0}^\delta+\big(I-M_0^{(K+1)M}\big)B^{-1} \zeta.
\end{align*}
This gives the desired expression of $\E[\hat x_{KM}^\delta]$.
Next, the variance component $\hat{e}_{(K+1)M}^\delta-\E[\hat{e}_{(K+1)M}^\delta]$ is given by
\begin{align*}
&\hat{e}_{(K+1)M}^\delta-\E[\hat{e}_{(K+1)M}^\delta]=M_0^M(\hat{e}_{KM}^\delta-\E[\hat{e}_{KM}^\delta])+c_0\sum_{i=0}^{M-1}H_{KM+i}M_0^i\hat{e}_{KM}^\delta\\
&+c_0\sum_{i=1}^{M}M_0^{M-i}(\zeta_{KM+i-1}-\zeta)+c_0^2\sum_{i=1}^{M-1}\sum_{t=0}^{M-i-1}H_{KM+i+t}M_0^t \zeta_{KM+i-1}\\
=&c_0\sum_{j=0}^K\sum_{i=0}^{M-1} M_0^{(K-j)M} H_{jM+i}M_0^i\hat{e}_{jM}^\delta
+c_0\sum_{j=0}^K\sum_{i=1}^{M} M_0^{(K-j+1)M-i}(\zeta_{jM+i-1}-\zeta)\\
&+c_0^2\sum_{j=0}^K\sum_{i=1}^{M-1}\sum_{t=0}^{M-i-1} M_0^{(K-j)M}H_{jM+i+t}M_0^t (\zeta_{jM+i-1}-\zeta)\\
&+c_0^2\sum_{j=0}^K\sum_{i=1}^{M-1}\sum_{t=0}^{M-i-1} M_0^{(K-j)M}H_{jM+i+t}M_0^t \zeta.
\end{align*}
Then it follows from the identity \eqref{eqn:sum-M} that
\begin{align*}
c_0\sum_{i=1}^{M-1}\sum_{t=0}^{M-i-1} H_{jM+i+t}M_0^t=&c_0\sum_{i=1}^{M-1}H_{jM+i}\Big(\sum_{t=0}^{i-1}M_0^t\Big)
=\sum_{i=1}^{M-1}H_{jM+i}(I-M_0^i) B^{-1}.
\end{align*}
Finally we derive
\begin{align*}
  &\hat{e}_{(K+1)M}^\delta-\E[\hat{e}_{(K+1)M}^\delta]\\
=&c_0\sum_{j=0}^K\sum_{i=0}^{M-1} M_0^{(K-j)M} H_{jM+i}M_0^i\hat{e}_{jM}^\delta
+c_0\sum_{j=0}^K\sum_{i=0}^{M-1} M_0^{(K-j+1)M-i-1}(\zeta_{jM+i}-\zeta)\\
&+c_0^2\sum_{j=0}^K\sum_{i=0}^{M-2}\sum_{t=0}^{M-i-2} M_0^{(K-j)M}H_{jM+i+t+1}M_0^t (\zeta_{jM+i}-\zeta)\\
&+c_0\sum_{j=0}^K \sum_{i=1}^{M-1}M_0^{(K-j)M}H_{jM+i}(I-M_0^i) B^{-1}\zeta\\
=&c_0\sum_{j=0}^K\sum_{i=0}^{M-1} M_0^{(K-j)M}\Big( H_{jM+i}\big(M_0^i\hat{e}_{jM}^\delta
+(I-M_0^i) B^{-1}\zeta\big)+ M_0^{M-i-1}(\zeta_{jM+i}-\zeta)\Big)\\
&+c_0^2\sum_{j=0}^K\sum_{i=0}^{M-2}\sum_{t=0}^{M-i-2} M_0^{(K-j)M}H_{jM+i+t+1}M_0^t (\zeta_{jM+i}-\zeta).
\end{align*}
This completes the proof of the proposition.

\subsection{Proof of Lemma \ref{lem:weierr_svrg}}

The proof employs the standard bias-variance decomposition and certain independence.
By Proposition \ref{prop:bias-var}, the following identities hold
\begin{align*}
&\E[R_1 (e_{(K+1)M}^\delta-B^{-1}\zeta)+R_2|\mathcal{F}^c_{(K+1)M}]=R_1(M_0^{(K+1)M}e_0^\delta-B^{-1} \zeta)+R_2,\\
&R_1 (e_{(K+1)M}^\delta-B^{-1}\zeta)+R_2-\E[R_1 (e_{(K+1)M}^\delta-B^{-1}\zeta)+R_2|\mathcal{F}^c_{(K+1)M}]\\
=&R_1(e_{(K+1)M}^\delta-\E[e_{(K+1)M}^\delta])=R_1\sum_{j=0}^K M_0^{(K-j)M} L_j(\zeta-B e_{jM}^\delta),
\end{align*}
where the random matrices $L_j$ are defined in \eqref{eqn:Lj}. Then we claim the following identity for any $i,i'=0,\ldots,M-1$,
\begin{equation}\label{eqn:svrg-ind}
  \E[\langle H_{jM+i}e_{jM}^\delta,H_{j'M+i'}e_{j'M}^\delta \rangle] = 0, \quad \mbox{if } i\neq i' \mbox{ or } j \neq j'.
\end{equation}
Clearly, it suffices to analyze the two cases $0\leq i<i'\leq M-1$, and $j<j'$ and $0\leq i,i'\leq M-1$
separately. Indeed, for any $0\leq i<i'\leq M-1$, the random matrix $N_{jM+i}$ is independent of
$G_{jM+i+1}e_{jM}^\delta$ and $N_{jM+i'}G_{jM+i'+1}e_{jM}^\delta$. Thus, using the identity
$\E_{jM+i}[N_{jM+i}]=0$, for any $i=0,\ldots, M-1$, we obtain
\begin{align*}
&\E_{jM+i}[\langle H_{jM+i}e_{jM}^\delta,H_{jM+i'}e_{jM}^\delta \rangle]\\
=&\E_{jM+i}[\langle N_{jM+i} G_{jM+i+1}e_{jM}^\delta,N_{jM+i'} G_{jM+i'+1}e_{jM}^\delta \rangle]\\
=& \langle \E_{jM+i}[N_{jM+i}] G_{jM+i+1}e_{jM}^\delta,N_{jM+i'} G_{jM+i'+1}e_{jM}^\delta \rangle=0. 
\end{align*}
Similarly, for any $j<j'$ and $0\leq i,i'\leq M-1$, the random matrix $N_{j'M+i'}$ is independent of
$N_{jM+i}G_{jM+i+1}e_{jM}^\delta$ and $G_{j'M+i'+1}e_{j'M}^\delta$, and hence
\begin{align*}
&\E_{j'M+i'}[\langle H_{jM+i}e_{jM}^\delta,H_{j'M+i'}e_{j'M}^\delta \rangle]\\
=&\E_{j'M+i'}[\langle N_{jM+i} G_{jM+i+1}e_{jM}^\delta,N_{j'M+i'} G_{j'M+i'+1}e_{j'M}^\delta \rangle]\\
=&\langle N_{jM+i} G_{jM+i+1}e_{jM}^\delta,\E_{j'M+i'}[N_{j'M+i'}] G_{j'M+i'+1}e_{j'M}^\delta \rangle=0. 
\end{align*}
The desired claim \eqref{eqn:svrg-ind} follows by taking full conditional of the
last two identities. Note that by assumption, $R_1$ is independent of $e_{(K+1)
M}^\delta-\E[e_{(K+1)M}^\delta]$. Then the bias-variance decomposition and the
claim \eqref{eqn:svrg-ind} imply
\begin{align*}
&\E[\E[\|R_1 (e_{(K+1)M}^\delta-B^{-1}\zeta)+R_2\|^2|\mathcal{F}^c_{(K+1)M}]]\\
=&{\rm I_0}+\E[\|R_1\sum_{j=0}^K M_0^{(K-j)M} L_j^1(\zeta-B e_{jM}^\delta)\|^2]\\
=&{\rm I_0}+c_0^2\sum_{j=0}^K\sum_{i=1}^{M-1}\E[\|R_1 M_0^{(K-j)M} H_{jM+i}(I-M_0^i)( e_{jM}^\delta-B^{-1}\zeta)\|^2].
\end{align*}
This and the definitions of the terms ${\rm I}_0$ and ${\rm I}_{1,j}$ complete the proof of the lemma.

\subsection{Proof of Lemma \ref{lem:weierr_sgd}}

The proof of the lemma is similar to Lemma \ref{lem:weierr_svrg}, {and employs suitable
independence relation crucially.} By Proposition \ref{prop:bias-var_sgd} and the standard bias-variance
decomposition, we have
\begin{align*}
&\E[\|R_1 (\hat{e}_{(K+1)M}^\delta-B^{-1}\zeta)+R_2\|^2]={\rm I_{0}}+
\E[\|R_1(\hat{e}_{(K+1)M}^\delta-\E[\hat{e}_{(K+1)M}^\delta])\|^2],
\end{align*}
with
$$\hat{e}_{(K+1)M}^\delta-\E[\hat{e}_{(K+1)M}^\delta]:=\sum_{j=0}^K \sum_{i=0}^{M-1} d_{j,i},$$
where $d_{j,i}$, in view of Proposition \ref{prop:bias-var_sgd}, are given by
\begin{align*}
d_{j,i}=&\mathrm{sgn}(M-1-i)c_0^2\sum_{t=0}^{M-i-2} M_0^{(K-j)M}H_{jM+i+t+1}M_0^t (\zeta_{jM+i}-\zeta)\\
&+c_0M_0^{(K-j)M}\big( H_{jM+i}\big(M_0^i\hat{e}_{jM}^\delta
+(I-M_0^i) B^{-1}\zeta\big)\\
 &+ M_0^{M-i-1}(\zeta_{jM+i}-\zeta)\big):=\sum_{t=0}^{M-i-2}d_{j,i,t}+d_{j,i,-1},
\end{align*}
where the notation $\text{sgn}(\cdot)$ denotes the sign function with the convention
$\text{sgn}(0)=0$. Next we repeat the argument for deriving \eqref{eqn:decom-svrg-1j},
and claim that $\E_{jM+i}[R_1 d_{j,i}]=0$ and $d_{j,i}|\mathcal{F}_{jM+i}\cup\mathcal{F}^c_{jM+i+1}$ is independent of
$d_{j',i'}|\mathcal{F}_{jM+i}\cup\mathcal{F}^c_{jM+i+1}$ for any $j\neq j'$ or $i\neq i'$ where $0\leq j'\leq j\leq K$,
$0\leq i,i'\leq M-1$. Indeed, the random vairable $d_{j',i'}$ is measurable with respect to $\mathcal{F}_{jM+i}\cup\mathcal{F}^c_{jM+i+1}$.
Then the direct computation using the identities $\E_{jM+i}[\zeta_{jM+i}-\zeta]=0$ and $\E_{jM+i}[H_{jM+i}]=0$ implies that
for any $0\leq j\leq K$ and $0\leq i\leq M-1$, the following identity holds
\begin{align*}
\E_{jM+i}[R_1 d_{j,i}]
=&\mathrm{sgn}(M-1-i)c_0^2\sum_{t=0}^{M-i-2}R_1 M_0^{(K-j)M}H_{jM+i+t+1}M_0^t \E_{jM+i}[\zeta_{jM+i}-\zeta]\\
&+c_0R_1 M_0^{(K-j)M}\big( \E_{jM+i}[H_{jM+i}]\big(M_0^i\hat{e}_{jM}^\delta
+(I-M_0^i) B^{-1}\zeta\big)\\
&+ M_0^{M-i-1}\E_{jM+i}[\zeta_{jM+i}-\zeta]\big)=0.
\end{align*}
Thus we derive
$$\E[\| R_1 (\hat{e}_{(K+1)M}^\delta-\E[\hat{e}_{(K+1)M}^\delta])\|^2]=\sum_{j=0}^K\sum_{i=0}^{M-1}\E[\| R_1 d_{j,i}\|^2].$$
Similarly, for fixed $j$, $i$ and any $0\leq t,t'\leq M-i-2$,  $\E[d_{j,i,t}|\mathcal{F}_{jM+i+t+1}]=0$ and $d_{j,i,t}|\mathcal{F}_{jM+i+t+1}$ is independent of $d_{j,i,t'}|\mathcal{F}_{jM+i+t+1}$ when $t> t'$. Consequently,
\begin{align*}
\E[\|R_1 d_{j,i}\|^2]=\sum_{t=-1}^{M-i-2}\E[\|R_1 d_{j,i,t}\|^2].
\end{align*}
Thus, we obtain
\begin{align*}
&\E[\|R_1(\hat{e}_{(K+1)M}^\delta-\E[\hat{e}_{(K+1)M}^\delta])\|^2]\\
=&c_0^2\sum_{j=0}^K\sum_{i=0}^{M-1}\E[\|R_1 M_0^{(K-j)M}\big( H_{jM+i}\big(M_0^i\hat{e}_{jM}^\delta
+(I-M_0^i) B^{-1}\zeta\big)+ M_0^{M-i-1}(\zeta_{jM+i}-\zeta)\big)\|^2]\\
&+c_0^4\sum_{j=0}^K\sum_{i=0}^{M-2}\sum_{t=0}^{M-i-2}\E[\|R_1 M_0^{(K-j)M}H_{jM+i+t+1}M_0^t (\zeta_{jM+i}-\zeta)\|^2].
\end{align*}
Reorganizing the last summation gives
\begin{align*}
&\sum_{i=0}^{M-2}\sum_{t=0}^{M-i-2}\E[\|R_1 M_0^{(K-j)M}H_{jM+i+t+1}M_0^t (\zeta_{jM+i}-\zeta)\|^2]\\
=&\sum_{i=1}^{M-1}\sum_{t=0}^{i-1}\E[\|R_1 M_0^{(K-j)M}H_{jM+i}M_0^t (\zeta_{jM+i-1-t}-\zeta)\|^2].
\end{align*}
This completes the proof of the lemma.


\bibliographystyle{abbrv}
\bibliography{sgd}
\end{document}